\documentclass[journal,twoside,web]{IEEEtran}
\usepackage[english]{babel}
\usepackage[utf8]{inputenc}
\usepackage{amsmath,amsfonts}
\usepackage{algorithmic}
\usepackage{array}
\usepackage[caption=false,font=normalsize,labelfont=sf,textfont=sf]{subfig}
\usepackage{textcomp}
\usepackage{stfloats}
\usepackage{url}
\usepackage{verbatim}
\usepackage{xcolor}
\usepackage{graphicx}  
\def\BibTeX{{\rm B\kern-.05em{\sc i\kern-.025em b}\kern-.08em
    T\kern-.1667em\lower.7ex\hbox{E}\kern-.125emX}}
\usepackage{balance}
\usepackage{multirow}
\usepackage{tikz}
\usepackage{listofitems}

\usepackage{makecell}

\usepackage{etoolbox}
\makeatletter
\@ifundefined{color@begingroup}%
  {\let\color@begingroup\relax
   \let\color@endgroup\relax}{}%
\def\fix@ieeecolor@hbox#1{%
  \hbox{\color@begingroup#1\color@endgroup}}
\patchcmd\@makecaption{\hbox}{\fix@ieeecolor@hbox}{}{\FAILED}
\patchcmd\@makecaption{\hbox}{\fix@ieeecolor@hbox}{}{\FAILED}

\def\BibTeX{{\rm B\kern-.05em{\sc i\kern-.025em b}\kern-.08em
    T\kern-.1667em\lower.7ex\hbox{E}\kern-.125emX}}

\usepackage{bm}

\newcommand{\herm}{{\scriptstyle \boldsymbol{\mathsf{H}}}}
\newcommand{\trans}{{\scriptstyle \boldsymbol{\mathsf{T}}}}

\usepackage{hyperref}

\newcommand{\C}{\mathbb C}

\newcommand{\CC}{\mathbf{c}}

\newcommand{\RR}{\mathbf{r}}
\newcommand{\XX}{\mathbf{x}}
\newcommand{\YY}{\mathbf{y}}
\newcommand{\PP}{\mathbf{p}}
\newcommand{\QQ}{\mathbf{q}}

\newcommand{\ZZ}{\mathbf{z}}

\newcommand{\Ad}{\mathbf A}

\newcommand{\Fd}{\mathbf F}

\newcommand{\Id}{\mathbf I}
\newcommand{\Cd}{\mathbf C}
\newcommand{\Ed}{\mathbf E}

\newcommand{\Kd}{\mathbf K}

\newcommand{\LLambda}{\boldsymbol{\Lambda}}

\usepackage{xcolor}
\definecolor{darkgreen}{rgb}{0.0, 0.5, 0.0}

\usetikzlibrary{arrows.meta,spy,calc}
\usepackage{xparse} 

\newcommand{\Thickline}{\\[-1.2em]\noalign{\hrule height 1.0pt}\\[-1.2em]}

\begin{document}

\title{
    DREAM: Deep-Reparametrization of Adaptive Regularization Maps for Fast Zero-Shot Self-Supervised Learning
}

\author{Thanh Trung Vu, Ander Biguri, Christoph Kolbitsch, Luca Calatroni, Kostas Papafitsoros,  Andreas Kofler
    \thanks{Corresponding author: Andreas Kofler {andreas.kofler[at]ptb.de} }%
    \thanks{Christoph Kolbitsch and Andreas Kofler are with Physikalisch-Technische Bundesanstalt (PTB), Braunschweig and Berlin, Germany.
    Thanh Trung Vu and Ander Biguri are with the University of Cambridge, Cambridge, UK,
    Luca Calatroni is with the University of Genoa, Genoa, Italy,
    Kostas Papafitsoros is with Queen Mary University of London, London, UK
    }%
}

\maketitle

\begin{abstract}
Adaptive regularization is an effective means of improving the flexibility of classical variational reconstruction methods while retaining their interpretability and mathematical structure. 
In this work, we propose an approach for 
learning
deep-reparameterized
adaptive regularization maps (DREAM)
for 
Total Variation (TV) and Total Generalized Variation (TGV) 
through algorithm unrolling in a zero-shot, self-supervised setting, requiring no access to large paired datasets. 
We first validate 
DREAM
on a two-dimensional image-denoising problem and then apply it to large-scale dynamic MRI reconstruction. The resulting adaptive models are theoretically grounded, as they retain the convergence properties of the primal–dual algorithm used as the building block of the unrolled architecture. Our results show that the regularization parameter maps can be learned efficiently, requiring substantially fewer weight updates than state-of-the-art methods to achieve comparable performance. We further demonstrate that the CNN parametrization of the maps acts as an implicit prior and provides effective dimensionality reduction, outperforming direct optimization of the adaptive parameters. Across both applications, the performance gap between supervised and zero-shot training is smaller for DREAM than for the end-to-end deep-learning methods considered for comparison. Remarkably, in dynamic cardiac MRI reconstruction, zero-shot training matches the performance of the supervised model. Moreover, the learned adaptive maps are highly interpretable, providing a compelling alternative to purely deep-learning-based methods when reference data are difficult to obtain.

\end{abstract}

\section{Introduction}

Over the last few years, learned reconstruction methods have emerged as state-of-the-art in the context of image reconstruction for various inverse problems \cite{ongie2020deep}, such as image denoising \cite{elad_image_2023}, computed tomography (CT) \cite{kulathilake2023review}, or Magnetic Resonance Imaging (MRI) \cite{chen2022ai}. 
There exist many different possibilities to employ neural networks for image reconstruction; a major distinction can be made based on whether the learning is physics-aware or not; see, for example, the recent review \cite{kofler2024quantitative} in the context of quantitative MRI.
Optimization-driven methods such as Plug-and-Play (PnP) approaches, which were pioneered in  \cite{venkatakrishnan2013plug}, or regularization by denoising (RED) approaches \cite{romano2017little}, can be used to solve arbitrary reconstruction problems by employing off-the-shelf Gaussian denoisers replacing first-order quantities (the gradient or the proximal operator) related to image regularization models. These denoisers can be parametrized in terms of large neural networks, such as the U-Net \cite{ronneberger2015u} or residual architectures like DnCNN \cite{zhang2017beyond}. This type of method is particularly attractive from a computational point of view, since the learning process can be entirely decoupled from the process of solving the actual inverse problem at hand.

Another particularly powerful class of learning-based methods is based on the so-called algorithm unrolling \cite{monga2021algorithm}, which was originally introduced for the task of improving and accelerating sparse coding solvers \cite{gregor2010learning}.
Unrolled methods interpret a finite (typically, small) number of iterations of an iterative algorithm as a neural network
whose parameters can be learned from data. This perspective allows retaining model/physics-based components (such as the forward model) modeling the acquisition of the observed data. 
Depending on the algorithmic building block used for their derivation, several reconstruction schemes based on, e.g.,\ ADMM \cite{sun2016deep}, Primal-Dual \cite{adler2018learned}, proximal gradient methods \cite{zhang2018ista}, etc. Nowadays,  unrolled models consistently occupy the first positions in organized image reconstruction challenges, see e.g.\ \cite{muckley2021results,beauferris2022multi,sidky2022report,lyu2025state}.
In terms of computational burdens,
the GPU memory necessary for end-to-end training linearly scales with the number of unrolled iterations and may pose limits to the usable network architectures. 
As an alternative, deep equilibrium models \cite{bai2019deep} were proposed in the context of inverse problems in \cite{gilton2021deep} to address this limitation by considering fixed-point mappings with weight tying, thus making the learning memory costs constant. 
The application of such approaches in the field of image reconstruction is vast \cite{Zou2023}, but still subject to significant training costs whenever fixed-point solvers have to be employed \cite{Daniele2026}.

A major limitation of neural network-based approaches is their lack of interpretability. For example, 
learned reconstruction methods by algorithmic unrolling are often hard to interpret as solvers of an underlying minimization problem without any further condition on the network architecture. In the context of MRI reconstruction, for instance, the iterates of the state-of-the-art Model-based Deep Learning (MoDL) approach developed in \cite{aggarwal2018modl} can be shown to correspond to the iterates of RED \cite{romano2017little} only whenever the employed neural network has a symmetric Jacobian, as pointed out in \cite{reehorst2018regularization}. Such conditions can be hard to verify in practice and may limit the expressivity of the proposed reconstruction scheme. Another major challenge of data-driven approaches is their need for large and diverse datasets for successful training. In a supervised training setting, datasets preferably consist of paired input-target data. However, for some applications, target data can be difficult or, in some cases, even impossible to obtain. 
To overcome such limitations, self-supervised training strategies have also been proposed in the past, see e.g.\
\cite{li_self-supervised_2025,wang_benchmarking_2026,tachella2026self} and references therein.
Many self-supervised methods for different image reconstruction problems have their root in image denoising (see, e.g.,  \cite{zhang_unleashing_nodate}). In the context of MRI reconstruction, the Self-Supervision by Data Undersampling (SSDU) framework \cite{yaman2020self} builds on the Noise2Self \cite{noise2self} denoising framework where training is carried out by randomly splitting the measurement data and using one set for reconstruction and the other set for the computation of a loss function, whose minimization drives the learning procedure. In \cite{yaman2020self} it was shown that when the access to \textit{many} different measurements is available, training neural network-based models using the SSDU framework is possible, and, in many cases, one can achieve nearly the same performance as for models trained in a supervised manner. However, although carried out in the absence of target data, training still requires a large number of different samples to succeed. To address this issue, zero-shot self-supervised approaches, where the model is trained solely on the data under consideration, were proposed \cite{yaman2022zeroshot}.
Learning expressive deep neural networks on a single sample requires dedicated strategies to avoid overfitting. For example, in the deep image prior  framework \cite{ulyanov2018deep}, the desired solution is reparametrized as a CNN and early stopping of the training algorithm is employed to enforce implicit regularization.  The deep image prior framework can be extended to include further hand-crafted regularization terms \cite{liu2019image}. Note that employing a simpler network architecture can also be used to avoid overfitting \cite{heckel2018deep}, otherwise, an overfitting mitigation mechanism strategy is required. In the Zero-Shot SSDU (ZS-SSDU) framework \cite{yaman2022zeroshot}, SSDU was extended by involving a third set of measurements to self-validate and provide a practical criterion for early stopping for training. 

\smallskip

In this work, we focus on learning an interpretable data-driven model from scratch, targeting situations in which target data are unavailable and/or large measurement datasets are difficult to obtain. In particular, we address the problem of learning adaptive regularization maps through algorithmic unrolling. We first validate DREAM on an illustrative two-dimensional image-denoising problem and then consider the challenging task of dynamic cardiac magnetic resonance imaging (MRI) reconstruction. In this setting, owing to the intrinsic constraints of the acquisition system and patient-specific limitations, such as limited breath-hold capability, data can generally be acquired only through accelerated scans. The resulting measurements are therefore undersampled and may contain motion and other acquisition artifacts, preventing the corresponding reconstructions from serving as reliable ground-truth references.
In this regard, even the reference images provided for supervised training in the recently released cardiac MRI datasets \cite{wang2024cmrxrecon,wang2025cmrxrecon2024} associated with the CMRxRecon challenge \cite{lyu2025state} are not reconstructed from fully sampled acquisitions. Rather, they are obtained from data acquired with acceleration factors of \(R=2\) and \(R=3\), respectively; see \cite{wang2024cmrxrecon,wang2025cmrxrecon2024} for details. Thus, although these images serve as training targets, they cannot be regarded as genuine fully sampled ground-truth references, motivating the need for self-supervised approaches.
We focus, in particular, on a hybrid approach that combines a hand-crafted regularization prior with a learnable parameter map that adapts the regularization strength to the local image features. This map is generated from the input image by a neural network and learned in an SSDU framework, without requiring ground-truth data. At test time, the trained network predicts an adaptive regularization map for each previously unseen input, which is subsequently used within an appropriate reconstruction algorithm to obtain the final image. In \cite{kofler2023learning}, such an idea was applied to learn an adaptive total-variation  (TV) regularization map \cite{pragliola2023and} in the context of dynamic cardiac MRI, quantitative MRI, dynamic image denoising, as well as 2D 
CT. 
In \cite{vu2025deep}, the method was extended to Total Generalized Variation (TGV) regularization \cite{BrediesKunischPock2010TGV} and applied to accelerated 2D MRI and 2D image denoising. 
Beyond gradient-based regularization, in \cite{kofler2025ell1}, a similar approach was applied to low-field MRI for learning adaptive sparsity level maps to be used for convolutional synthesis regularization, while in \cite{schulz2026learning}, 
the method was demonstrated to be more robust with respect to data-distribution shifts. 
Importantly, in all these methods, the networks were trained in a supervised setting, partly with considerable training times of several days, as well as noticeable hardware requirements.

Here, differently from \cite{kofler2023learning} and \cite{vu2025deep}, we consider the problem of 
adaptive learning TV and TGV regularization maps 
parametrized 
by 
a U-Net \cite{ronneberger2015u} 
in a zero-shot self-supervised setting.
By some simple but efficient adaptations, we demonstrate that the combination of adaptive TV and TGV priors enables fast and efficient zero-shot learning using the SSDU paradigm. Crucially, we show that for the considered dynamic MRI example, the obtainable results are in fact on par with the ones obtained by supervised training as well as supervised training with additional test-time fine-tuning. This importantly implies that no training data at all is required for the investigated approaches. On the other hand, these approaches could be used to generate reconstructions that could serve as target images for other supervised training methods. Last,
the proposed approach, which we call DREAM, 
does not require early stopping or particular network constraints since the role of the network employed is, by design, limited to providing the adaptive regularization strength for the regularizers considered.

\begin{figure}
    \centering
    \includegraphics[width=\linewidth]{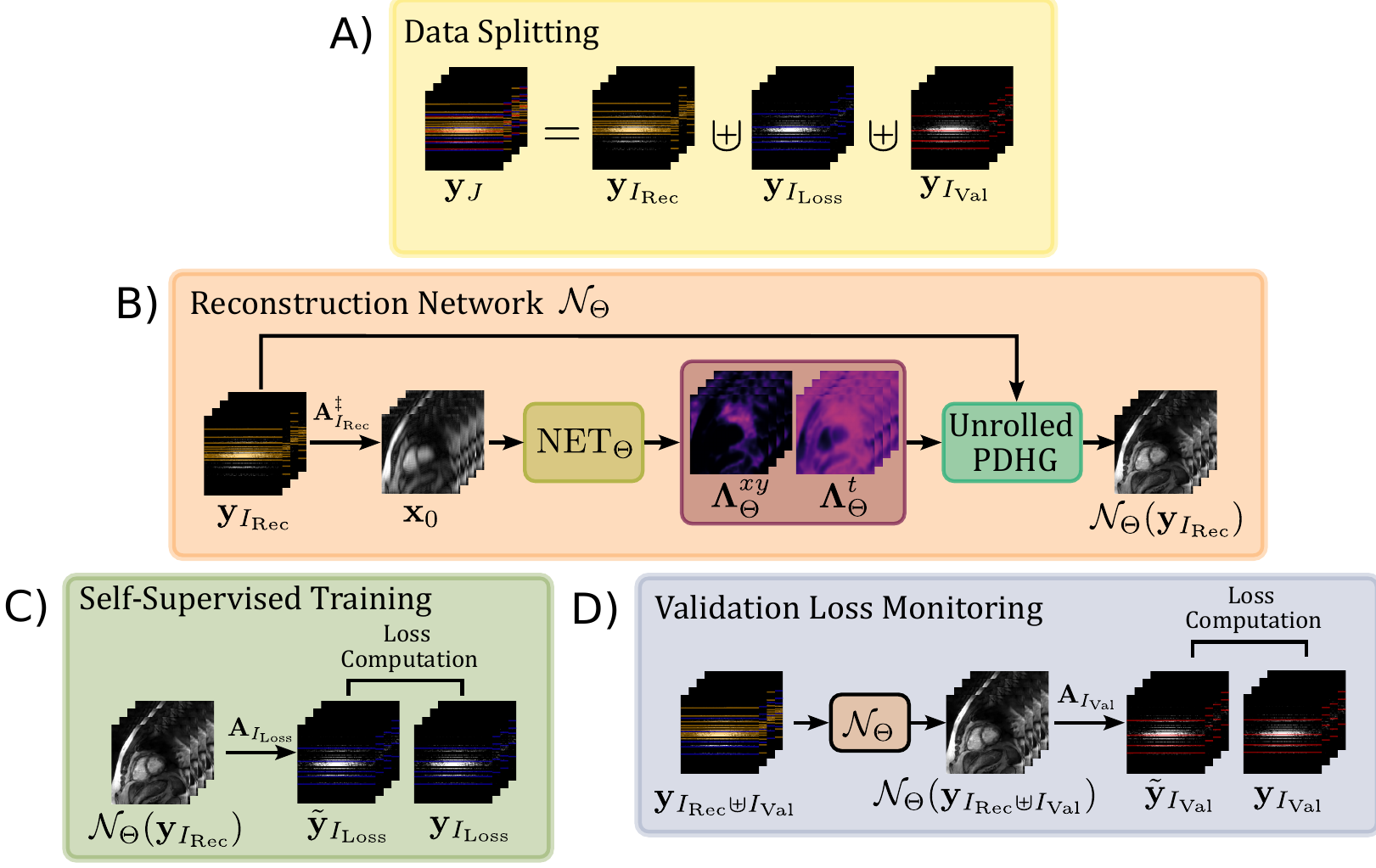}
    \caption{Proposed Zero-Shot Test-Time Training (ZS-TTT) approach for learning adaptive regularization parameter maps (only shown for TV-$\LLambda$ for a dynamic cardiac MRI problem). A) The raw measurement data $\YY_J$ is split into three disjoint sets $\YY_{I_{\mathrm{Rec}}}$, $\YY_{I_{\mathrm{Loss}}}$ and $\YY_{I_{\mathrm{Val}}}$, which are used for reconstruction, loss computation and validation. B) The data $\YY_{I_{\mathrm{Rec}}}$ is first reconstructed to an initial image by a reconstruction operator $\Ad_{I_{\mathrm{Rec}}}^\ddagger$, from which the network $\mathrm{NET}_{\Theta}$ estimates adaptive regularization parameter maps (here $\LLambda_{\Theta}:=[\LLambda_{\Theta}^{xy},\LLambda_{\Theta}^{xy},\LLambda_{\Theta}^{t}]^\trans$). These are used for the definition of a problem formulation (e.g.\ \eqref{eq:tv_problem_Lambda} or \eqref{eq:tgv_problem_Lambda}), which is then solved by an unrolled primal-dual hybrid gradient algorithm (PDHG). When no target image is available, the reconstructed image can be transformed to the measurement space by the operator $\Ad_{I_{\mathrm{Loss}}}$. C) The estimated measurement data $\tilde{\YY}_{I_{\mathrm{Loss}}}$ are compared to all data \textit{not} used for reconstruction, and a loss function is computed as in \eqref{eq:loss_function_ssdu}. D) Finally, if necessary, overfitting can be avoided by monitoring the difference (e.g.\ in terms of the MSE) between the estimated measurement data restricted to the validation set $\Ad_{I_{\mathrm{Val}}} \mathcal{N}_{\Theta}(\YY_{I_{\mathrm{Rec}}\uplus I_{\mathrm{Val}}})$ and the validation data $\YY_{I_{\mathrm{Val}}}$. }

    \label{fig:approach}
\end{figure}

\section{Problem Formulation and Reconstruction Methods}\label{sec:problem_formulations_and_rec_methods}
We define in the following the reconstruction models considered. We focus in particular locally adaptive TV and TGV regularization models given their broad popularity in the context of image reconstruction problems, although a similar treatment may also be considered for different regularization models, see, e.g., \cite{ChambolleLions1997,holler2014infimal,lefkiammiatis2013} for higher-order regularization models and \cite{pourya2025dealing} for learning regularization maps in the case of Fields-Of-Experts-type regularizers \cite{Roth2005}.

\subsection{Notation}\label{subsec:notation}
Let $V,\, W \in \{\mathbb{R}, \mathbb{C}\}$ and $\Ad_J: V^N \rightarrow W^M, \XX \mapsto \YY$ denote a linear operator from the space of images to their measurement space, depending on $J$, which refers to some modality-specific set of parameters.  The forward problem we aim to solve reads:
\begin{equation}\label{eq:forward_problem}
  \text{find }\XX_{\mathrm{true}}\quad\text{s.t. }\quad \YY := \Ad_J \XX_{\mathrm{true}} + \mathbf{e},
\end{equation}
where $\mathbf{e}\sim \mathcal{N}(0,\sigma^2\mathbf{I})$ denotes an additive white Gaussian noise (AWGN) perturbation modeling measurement errors. Denoting by $\mathcal{R}_{\lambda}(\XX)$ a regularization function whose contribution is \textit{globally} weighted by a scalar parameter $\lambda>0$. and by $\mathcal{R}_{\LLambda}(\XX)$  a regularizer whose contribution is \textit{locally} adapted by means of a regularization parameter map $\LLambda$, we consider the two instances of the reconstruction problem:
\begin{align}
    \underset{\XX}{\min}\, \, \frac{1}{2}\| \Ad_J \XX - \YY \|_{2}^2 + \mathcal{R}_{\lambda}(\XX),  \tag{$P_{\lambda}$} \label{eq:variational_problem_lambda} \\
    \underset{\XX}{\min}\, \, \frac{1}{2}\| \Ad_J \XX - \YY \|_{2}^2 + \mathcal{R}_{\LLambda}(\XX),  \tag{$P_{\LLambda}$} \label{eq:variational_problem_Lambda}
\end{align}
where the quadratic fidelity models the presence of AWGN.

\subsection{Dynamic Cardiac MRI}\label{subsec:dynamic_cardiac_mri}
For image denoising, the forward operator in \eqref{eq:forward_problem} is simply the identity, while for dynamic cardiac MRI, it is more involved. In dynamic cardiac MRI, the sought image depicts the temporal evolution of a moving heart and is used for cardiac functional assessment, see e.g.\ \cite{kolbitsch2014cardiac}. The challenge in the reconstruction problem lies in its relatively high dimensionality from a computational point of view, on the one hand, and, on the other hand, in being able to faithfully depict the cardiac movement, which is the clinical feature of interest.
Let $\XX_{\mathrm{true}} \in \C^{N_{\mathrm{time}}\cdot N} $ with $N = N_x \cdot N_y$ denote the vector representation of the complex-valued dynamic target image of spatial dimension $N_x \times N_y$. The accelerated measurement process can be modeled by \eqref{eq:forward_problem}, where $\YY \in \C^{N_{\mathrm{coils}} \cdot N_{\mathrm{time}} \cdot M}$ denotes the dynamic multi-coil $k$-space data with $M<N$. The operator $\Ad_J$ denotes the forward operator that samples each 2D image $\XX_i, i=1,\ldots, N_{\mathrm{time}}$ in its 2D Fourier domain, and takes the form
\begin{equation}
    \begin{aligned}\label{eq:multi_coil_mri_fwd}
        \Ad_J: \C^{N_{\mathrm{time}}\cdot N}           & \rightarrow \C^{N_{\mathrm{coils}} \cdot N_{\mathrm{time}}\cdot M}, \\
        [\XX_1,\ldots, \XX_{N_{\mathrm{time}}}]^\trans & \longmapsto \bigoplus_{i=1}^{N_{\mathrm{time}}} \Ad_{J_i} \XX_i.
    \end{aligned}
\end{equation}
Thereby, each $\Ad_{J_i}$ given as
\begin{equation}
    \Ad_{J_i}:=(\Id_{N_{\mathrm{coils}} }\otimes \Ed_{J_i}) \Cd,
\end{equation}
where $\Ed_{J_i}:= \mathbf{U}_{J_i} \Fd$ denotes a 2D Fourier transform $\Fd$ followed by an undersampling mask $\mathbf{U}_{J_i}$ that only retains the Fourier coefficients index by the set $J_i$. Moreover, $\Cd:=[\Cd_1, \ldots, \Cd_{N_{\mathrm{coils}}}]^\trans$ with $\Cd_k:=\mathrm{diag}(\CC_k), \CC_k\in\C^N$ represents the so-called coil-sensitivity map (CSM) operator that weights the image by multiplying it with the respective coil-sensitivity profiles. In applications, the operator $\Cd$ is unknown at first and must be first estimated with adequate methods, e.g., \cite{walsh_adaptive_2000}, \cite{uecker2013}. After having computed an estimate $\tilde{\Cd}$ of $\Cd$ (and thus an estimate $\tilde{\Ad}_J$ of $\Ad_J$), one can proceed with the reconstruction.
Note that, typically, the set of coefficients $I_i$ acquired for each time point varies to collect complementary information and achieve incoherent sampling, which yields high-dimensional noise-like artefacts that can be eliminated by Compressed Sensing-like approaches, such as TV, TGV, and others \cite{lustig2007sparse}.

\subsection{Regularization Models}
\subsubsection{Adaptive TV}

The locally adaptive version of the TV reconstruction problem employing the locally adaptive regularization parameter map $\LLambda$ is obtained by setting $\mathcal{R}_{\LLambda}(\XX):=\| \LLambda \nabla \XX\|_1$ as opposed to $\mathcal{R}_{\lambda}(\XX):= \lambda \| \nabla \XX\|_1$ for $\lambda >0$. The resulting reconstruction problem \eqref{eq:variational_problem_Lambda} specifies into
\begin{equation}
    \underset{\XX}{\min}\, \, \frac{1}{2}\|\Ad_J \XX -  \YY\|_2^2 + \| \LLambda \nabla \XX\|_1  \tag{$\mathrm{TV}_{\LLambda}$}, \label{eq:tv_problem_Lambda}
\end{equation}
where, in our setting, $\|\cdot\|_1$ denotes the $\ell_1$-norm, so that the proposed regularizers rely on anisotropic versions of TV.

\subsubsection{Adaptive TGV}
In analogy to \eqref{eq:tv_problem_Lambda}, we consider a locally adaptive TGV model.
In the classical context, TGV generalizes TV by incorporating higher-order derivatives,
which helps preserve smooth transitions while still promoting 
edges
\cite{BrediesKunischPock2010TGV}.
In practice, the second-order TGV is most commonly used and, in its adaptive form, involves two regularization  maps
$\LLambda := \bigl[\LLambda_{0}, \LLambda_{1}\bigr]^\trans$,
replacing the scalar parameters $\lambda_0$ and $\lambda_1$ by locally adaptive weights, see \cite{vu2025deep}.
The corresponding anisotropic adaptive TGV model reads
\begin{equation}\label{eq:tgv_reg_Lambda}
\begin{aligned}
\mathcal{R}_{\LLambda}(\XX)
&:= \operatorname{TGV}_{\LLambda_{0},\,\LLambda_{1}}(\XX) \\
&:= \min_{\mathbf{w}}\;
\bigl\| \LLambda_{1}\,(\nabla \XX - \mathbf{w})\bigr\|_1
+ \bigl\| \LLambda_{0}\,\boldsymbol{\mathcal{E}}\mathbf{w}\bigr\|_1 ,
\end{aligned}
\end{equation}
where \(\mathbf{w}\)
is an auxiliary vector field and \(\boldsymbol{\mathcal{E}}\) denotes the (discrete) symmetrized gradient.
Using again the (squared) Euclidean distance as fidelity term, the resulting TGV reconstruction problem reads
\begin{equation}
\small{
\underset{\XX,\mathbf{w}}{\min}\, \,
\frac{1}{2}\|\Ad_J \XX -  \YY\|_2^2
+ \bigl\| \LLambda_{1}\,(\nabla \XX - \mathbf{w})\bigr\|_1
+ \bigl\| \LLambda_{0}\,\boldsymbol{\mathcal{E}}\mathbf{w}\bigr\|_1.
\tag{$\mathrm{TGV}_{\LLambda}$}
\label{eq:tgv_problem_Lambda}
}
\end{equation}

\subsection{Learned Reconstruction by Algorithmic Unrolling}
As a building block of the proposed learned procedure, we consider
the primal-dual hybrid gradient (PDHG) algorithm \cite{ChambollePock2011FirstOrderPrimalDual}, which is a well-known algorithm for solving problems of the form
\begin{equation}\label{eq:min_fkx_gx}
    \underset{\XX}{\min}\, f(\mathbf{K}(\XX)) + g(\XX),
\end{equation}
where $f,g$ are proper, convex, lower-semicontinuous functionals and $\mathbf{K}$ is a linear operator. 
The PDHG iteerations read:
\begin{equation}\label{eq:pdhg_iterates}
    \begin{aligned}
        \ZZ_{k+1}       & = & \mathrm{prox}_{\sigma\, f^\ast}(\ZZ_k + \sigma \Kd \bar{\XX}_k ), \\
        \XX_{k+1}       & = & \mathrm{prox}_{\tau\, g}(\XX_k - \tau \Kd^\herm \ZZ_{k+1} ),      \\
        \bar{\XX}_{k+1} & = & \XX_{k+1} + \theta(\XX_{k+1} - \XX_k),
    \end{aligned}
\end{equation}
where the step-sizes $\sigma, \tau$ fulfill $\sigma \tau < 1/\|\Kd\|_2^2$ and $\theta=1$, $f^\ast$ denotes the convex conjugate of $f$ and $\mathrm{prox}_{\sigma\, f^\ast}$ and $\mathrm{prox}_{\tau\, g}$ denote the proximal operators of the scaled functions $\sigma f^\ast$ and $\tau g$, respectively.

\subsubsection{PDHG for adaptive TV}
To match \eqref{eq:tv_problem_Lambda} with \eqref{eq:min_fkx_gx}, we set $f:=f_1+f_2$ given by
\begin{equation}\label{eq:matching_tv}
    f(\PP,\QQ):=f_1(\PP) + f_2(\QQ) = \frac{1}{2}\| \PP - \YY\|_2^2 + \|\LLambda \QQ\|_1,
\end{equation}
with the choice $\Kd:=[\Ad\, ,\, \nabla]^\trans$ and $g(\XX)\equiv 0$, see also \cite{sidky2012convex, kofler2023learning} for additional details.
Note that one needs to compute the proximal operator of the convex conjugate of the weighted $\ell_1$-norm with weights defined by $\LLambda$; this amounts to the application of a component-wise projection operator that projects the $i$-th component of the input onto the bilateral set $[-(\LLambda)_i, (\LLambda)_i]$ defined by the values in $\LLambda$, as detailed in \cite{kofler2023learning}.

\subsubsection{PDHG for TGV}
To match \eqref{eq:tgv_problem_Lambda} with \eqref{eq:min_fkx_gx}, let \(\mathbf{z}:=[\XX,\mathbf{w}]^\trans\) and define the stacked linear operator
\begin{equation}
\mathbf{K}(\XX,\mathbf{w})
:=
\bigl[\Ad_J \XX,\ \nabla \XX-\mathbf{w},\ \boldsymbol{\mathcal{E}}\mathbf{w}\bigr]^\trans.
\label{eq:tgv_K3_operator}
\end{equation}
We then set
$
g(\XX,\mathbf{w}) := 0,
$ and
\[
f(\PP,\QQ,\RR)
:=
\frac{1}{2}\|\PP-\YY\|_2^2
+\|\LLambda_{1}\QQ\|_1
+\|\LLambda_{0}\RR\|_1,
\label{eq:tgv_split_fKg_3layer}
\]
to have the problem 
\[
\min_{\XX,\mathbf{w}}\; f\bigl(\mathbf{K}(\XX,\mathbf{w})\bigr) + g(\XX,\mathbf{w}),
\]
which is in the form of \eqref{eq:min_fkx_gx}.
Compared to the TV case, the only additional ingredient is the presence of two weighted $\ell_1$ terms; accordingly, the dual updates involve two component-wise projections, with bounds given pointwise by $\LLambda_{1}$ and $\LLambda_{0}$.

\subsubsection{Algorithm Unrolling}
We follow the notation in  \cite{kofler2023learning} and denote by $\XX_{\LLambda}^T:= S^T(\XX_0, \YY, \Ad_J, \LLambda)$ the $T$-th iterate of PDHG for solving \eqref{eq:variational_problem_Lambda} for fixed $\LLambda$ with starting value $\XX_0$. Note that $\XX_{\LLambda}^T \rightarrow S^\ast(\YY, \Ad_J, \LLambda)$, where $S^{\ast}(\YY, \Ad_J, \LLambda)$ solves  \eqref{eq:variational_problem_Lambda}.\\
As in \cite{kofler2023learning}, we reparametrize the adaptive regularization parameter map $\LLambda$ as the output of a CNN applied to some input image $\XX_{\mathrm{input}}$, i.e., upon vectorization, we have $\LLambda:=\LLambda_{\Theta}:= \mathrm{NET}_{\Theta}(\XX_{\mathrm{input}})$, with $\Theta \in \mathbb{R}^{\ell}$ denoting the set of trainable parameters.
By $\mathcal{N}_{\Theta}^T$, we denote, by slight abuse of notation, the composition of the network $\mathrm{NET}_{\Theta}$ and the unrolling of $T>0$ iterations of the scheme $S^T(\XX_0, \YY, \Ad_J, \LLambda_{\Theta})$ to approximately solve \eqref{eq:variational_problem_Lambda}, i.e.\ $\mathcal{N}_{\Theta}^T(\YY):= S^T(\XX_0, \YY, \Ad_J, \mathrm{NET}_{\Theta}(\XX_{\mathrm{input}})\big)$, where $\XX_{\mathrm{input}}:=\Ad_J^\ddagger \YY$ is obtained by a reconstruction operator $\Ad_J^\ddagger$, e.g., $\Ad_J^\ddagger \in \{\Ad_J^\herm, \Ad_J^\trans, \Ad_J^\dagger, \ldots\}$ with $\Ad_J^\dagger$ being the Moore-Penrose pseudo-inverse.
Figure~\ref{fig:approach}B visualizes the reconstruction network $\mathcal{N}_{\Theta}^T$ for the dynamic cardiac MRI problem, with some notations to be defined in later sections.

\subsection{Self-Supervised Training via Data Undersampling (SSDU)}
Regularization parameter maps parametrized by CNNs have been successfully learned for TV \cite{kofler2023learning}, TGV \cite{vu2025deep}, as well as convolutional dictionaries \cite{kofler2025ell1}. In all these works, supervised training was adopted, i.e., 
assuming access to a dataset $\mathcal{D}$ of input-target data pairs given in the form \eqref{eq:forward_problem},
a suitable set of trainable parameters $\Theta^{\ast}$ was obtained by minimizing a loss function
\begin{equation}\label{eq:loss_function_sup}
    \mathcal{L}^{\mathrm{sup}}(\Theta) := \sum_{(\YY,\XX_{\mathrm{true}}) \in \mathcal{D}} l \big( \mathcal{N}_{\Theta}(\YY), \XX_{\mathrm{true}}\big) + r(\Theta).
\end{equation}
Here, $l(\,\cdot\,,\,\cdot\,)$ denotes some appropriate discrepancy or similarity metric between images (e.g.\ the mean squared error (MSE) loss, the structural similarity index measure (SSIM), etc.) and $r(\,\cdot\,)$ some regularizational functional for the weights, such as weight-decay, which can be used to prevent overfitting.

Here, we explore the possibility of training the network parameters $\Theta$ by self-supervision in a zero-shot fashion.
We focus on SSDU \cite{yaman2020self} and its zero-shot variant (ZS-SSDU) \cite{yaman2022zeroshot}, but other techniques such as Noise2Void \cite{krull2019noise2void}, Noise2Self \cite{noise2self} could be used as well.
As noted in \cite{yaman2020self}, although ZS-SSDU can be applied for subject-specific learning, it typically requires early stopping to avoid overfitting \cite{hosseini2020high}, similar to other approaches such as deep image prior.

Let $J$ denote the set of samples that the operator $\Ad_J$ acquires, e.g., the set of $k$-space coefficients in MRI. Further, let  $I_{\mathrm{Val}}, I_{\mathrm{Rec}},I_{\mathrm{Loss}} \subset J$ with $I_{\mathrm{Loss}}:= (J\setminus I_{\mathrm{Val}}) \setminus I_{\mathrm{Rec}}$. By this, we have the partition $I_{\mathrm{Rec}}\, \uplus  I_{\mathrm{Loss}}\, \uplus I_{\mathrm{Val}}=J$. 
A visualization of this data splitting applied to the acquired dynamic cardiac MRI measurement is shown in Figure~\ref{fig:approach}A.
In the SSDU framework, where no target images are required, the discrepancy term $l(\,\cdot\, , \,\cdot\,)$ in \eqref{eq:loss_function_sup} operates in the raw-data space and takes the form
\begin{equation}\label{eq:loss_term_ssdu}
    l\big(\Ad_{I_{\mathrm{Loss}}} \mathcal{N}_{\Theta}(\YY_{I_{\mathrm{Rec}}}), \YY_{I_{\mathrm{Loss}}}\big),
\end{equation}
where $\YY_{I_{\mathrm{Rec}}}$ and $\YY_{I_{\mathrm{Loss}}}$ denote the sets of measurements restricted to the sets $I_{\mathrm{Rec}}$ and $I_{\mathrm{Loss}}$, respectively.
Figure~\ref{fig:approach}C visualizes the components of this loss term applied to the MRI problem.
Thus, assuming only access to one set of measurements $\YY$, the overall zero-shot SSDU loss function is given by
\begin{equation}\label{eq:loss_function_ssdu}
    \mathcal{L}_{\mathrm{ZS}}^{\mathrm{SSDU}}(\Theta) := \sum_{( I_{\mathrm{Rec}}, {I_{\mathrm{Loss}}})} l\big(\Ad_{I_{\mathrm{Loss}}} \mathcal{N}_{\Theta}(\YY_{I_{\mathrm{Rec}}}), \YY_{I_{\mathrm{Loss}}}\big) + r(\Theta),
\end{equation}
where $\mathcal{S}_{J}$ denotes the set of possible partitions
of $J$, i.e., $\mathcal{S}_{J}:=\{ (I_{\mathrm{Rec}}, I_{\mathrm{Loss}}) \,|\,  I_{\mathrm{Rec}}, I_{\mathrm{Loss}} \subset J, \, I_{\mathrm{Rec}} \uplus I_{\mathrm{Loss}} \uplus I_{\mathrm{Val}} = J\}$.
Because the number of possible partitions grows exponentially with the cardinality of $J$, \eqref{eq:loss_function_ssdu} is infeasible to compute in practice. Therefore, choosing a subset of $\mathcal{S}_{J}$ is necessary in practice.
Note that, when multiple measurements data $\YY_1, \ldots, \YY_{N_{\mathrm{train}}}$ are available, \eqref{eq:loss_function_sup} involves a further summation over the different measurements. The minimization of \eqref{eq:loss_function_ssdu} is tackled similarly to supervised loss functions, i.e., by stochastic gradient descent techniques. Thereby, to avoid overfitting, the set $I_{\mathrm{Val}}$ serves as a validation set to monitor the quality of the reconstruction by comparing its estimated measurements to the ones restricted to $I_{\mathrm{Val}}$, i.e. by monitoring $l\big(\Ad_{I_{\mathrm{Val}}} \mathcal{N}_{\Theta}(\YY_{I_{\mathrm{Rec}}\uplus I_{\mathrm{Loss}}}), \YY_{I_{\mathrm{Val}}}\big)$. 
Figure~\ref{fig:approach}D illustrates this validation loss monitoring process.
After having obtained a suitable set of parameters $\Theta$, one can compute an estimate of the solution by applying the trained network to the entire given measurements.

\section{Experiments}

Here, we apply 
DREAM 
for both the TV-$\LLambda$ and TGV-$\LLambda$ regularizers to an illutrative image denoising example and a challenging dynamic cardiac MRI problem. 
All implementations use PyTorch, and the full code for all the experiments can be found at \url{https://github.com/trung-vt/ZS_SSL_UTGV}. For the MRI example, all operators (e.g.\ the Fourier, coil sensitivity operator, the finite differences operator $\nabla$, the symmetrized gradient $\boldsymbol{\mathcal{E}}$ etc.), as well as the implementations of the PDHG algorithm are taken from \texttt{MRpro} \cite{zimmermann2025mrpro}, \cite{zimmermann2026mrpro_arxiv}.

\subsection{Experimental Setup Design}
We will compare all methods in terms of their performance for the following set-ups:
\begin{enumerate}
    \item Supervised training;
    \item Supervised training + test-time training (TTT) for self-supervised fine-tuning using the SSDU loss in \eqref{eq:loss_function_ssdu}, similar to in \cite{hosseini2020high};
    \item Zero-Shot self-supervised training using SSDU. 
\end{enumerate}
Supervised training serves to evaluate the baseline performance that can be achieved given a suitable paired dataset for training. The second approach serves to adapt/fine-tune the pre-trained model to each unseen sample in the test data by further adapting the model using self-supervision. Note that this test-time fine-tuning to each sample of the test is particularly beneficial when a model is exposed to a data distribution shift \cite{darestani2022test}, which we, however, do not have in our setting. The third ZS-TTT experiment reveals the achievable performance of a model without having ever been exposed to training data.
For the self-supervised fine-tuning and zero-shot experiments, we further distinguish three cases. For the first case, we use a predetermined number of training iterations as a stopping criterion, representing the case where we keep training until reaching the practical computational limit. For the second case, we use an estimate of the validation SSDU loss, which is calculated on the validation set $I_{\mathrm{Val}}$ as a model selection criterion. For the third, the best model is selected based on the lowest MSE with respect to the unknown target image. The three different stopping strategies serve to investigate whether overfitting to the noise occurs for any model and, if yes, how much the second model selection criterion agrees with the third, which is ideal but not applicable in practice.

For image quality assessment, we use the peak signal-to-noise ratio (PSNR), the mean squared error (MSE), the structural similarity index measure (SSIM)  \cite{wang2004image}, and the Haar Wavelet-based perceptual similarity index (HaarPSI) \cite{reisenhofer2018haar}.

\subsection{Image Denoising}

\paragraph{Dataset}
We evaluate the denoising setting using a subset of 50 images from the DIV2K training dataset \cite{agustsson2017ntire}.
The images are processed as single-channel images with pixel intensities normalized to the range $[0,1]$ and cropped to size $512 \times 512$. Additive Gaussian 
noise with $\sigma \in \{0.05, 0.10, 0.15\}$ was used.
For the self-supervised loss, a split of $81\%/9\%/10\%$ of the pixels is used to define the sets $I_{\mathrm{Rec}}$, $I_{\mathrm{Loss}}$, and $I_{\mathrm{Val}}$.

\paragraph{Methods of Comparison}

As a method of comparison, we use the pre-trained ``Reconstruct Anything Model`` (RAM) denoiser \cite{terris2026ram}, which is provided by the DeepInv framework \cite{tachella2025deepinverse}. Although it is a general-purpose foundation model, it is still a powerful Gaussian denoiser and hence is chosen as a baseline for comparison. The model employs 35,618,813 trainable parameters.

\paragraph{Implementation Details}
The pretrained supervised adaptive TV-$\LLambda$ and TGV-$\LLambda$ models are the same as those tested in \cite{vu2025deep}.
For zero-shot training, all settings for adaptive TV and TGV models are also kept the same as in \cite{vu2025deep}, except the U-Net architecture \cite{ronneberger2015u} now being much smaller (more than 260 times fewer parameters), 
with two encoding and two decoding blocks, 16 initial feature channels,
resulting in an 1-16-32-64-32-16-$K$ UNet structure, where $K$ is the number of output channels ($K = 1$ for TV-$\LLambda$ and $K = 2$ for TGV-$\LLambda$). This means that for TV-$\LLambda$, the estimated parameter map takes the form $\LLambda_{\Theta}:=[\LLambda_{\Theta}^{xy}, \LLambda_{\Theta}^{xy}]^\trans$ to equally weight the spatial derivatives in the $x$- and $y$-directions. For TGV-$\LLambda$, the estimated parameter maps take the form 
$\LLambda_{1,\Theta}:=[\LLambda_{1,\Theta}^{xy}, \LLambda_{1,\Theta}^{xy}]^\trans$ and $\LLambda_{0,\Theta}:=[\LLambda_{0,\Theta}^{xy}, \LLambda_{0,\Theta}^{xy}]^\trans$. 
The respective networks $\mathrm{NET}_{\Theta}$ consist of 109,073 and 109,090 trainable parameters, respectively.

\paragraph{Training}

For all models, ZS-TTT is carried out for a predetermined maximum of $500$ epochs, representing the practical computational limit for one denoising example. 
TV-$\LLambda$ and TGV-$\LLambda$ were trained using Adam \cite{kingmaAdamMethodStochastic2015a} with a learning rate of $10^{-4}$ and no weight decay, while RAM is trained using AdamW \cite{loshchilov2018decoupled} with a learning rate of $2 \times 10^{-4}$ and weight decay $10^{-5}$.
The number of unrolled iterations for TV-$\LLambda$ and TGV-$\LLambda$ is $T=96$ and $T=128$, respectively.

\subsection{Dynamic Cardiac MRI}

\paragraph{Dataset}

We used the  MDCNN dataset \cite{mdcnn_data}, which consists of radially acquired $k$-space data measurements from 108 subjects (101 patients and 7 healthy subjects) obtained with a balanced steady state free precession (bSSFP sequence) on a Siemens Magnetom Vida scanner.  
The image dimensions are $N_x \times N_y \times N_t = 208 \times 208 \times 25$, and the number of CSMs ranges from 15 to 22; for further details see \cite{mdcnn_data}. To define target images to be used for supervised training and image quality assessment, the full $k$-space data was reconstructed with an iterative SENSE approach \cite{pruessmann2001advances}, 
by estimating the CSMs from the given $k$-space measurements using the method \cite{walsh_adaptive_2000}. 
From these images, we retrospectively generated undersampled multi-coil Cartesian dynamic $k$-space data by sub-sampling the phase encoding direction using a Gaussian variable density sampling pattern for three acceleration factors \(R_\text{train} \in \{4, 6, 8\}\). Thereby, the $n_{\mathrm{center}}=10$ lines were always sampled for each temporal frame, and the $k$-space data was corrupted by complex-valued Gaussian noise with a noise level \(\sigma_{\mathrm{train}} = 0.05\). The CSMs to be used within the reconstruction were estimated from the 24 simulated calibration $k$-space lines. By doing so, the CSMs used for the simulation and reconstruction differ from each other, representing a realistic model-mismatch that is encountered in practice.

\paragraph{Methods of Comparison}
As a representative method of comparison, we choose the well-established Model-based Deep Learning (MoDL) approach \cite{aggarwal2018modl}, which is a physics-informed method that unrolls $T$ iterations of a scheme for which the application of the network module can be interpreted as a learned proximal operator used within a half-quadratic splitting scheme \cite{kofler2024quantitative}. Our implementation of MoDL has 338,689 trainable parameters.

\paragraph{Implementation Details}
TV-$\LLambda$ uses a 3D U-Net architecture \cite{ronneberger2015u} with
two encoding and two decoding blocks, 16 initial filters, and 2 output channels (one for the spatial parameter map, one for the temporal parameter map),
resulting in a 2-16-32-64-32-16-2 structure. This means that for TV-$\LLambda$, the estimated parameter map takes the form $\LLambda_{\Theta}:=[\LLambda_{\Theta}^{t}, \LLambda_{\Theta}^{xy}, \LLambda_{\Theta}^{xy}]^\trans$ to equally weight the spatial derivatives in the $x$- and $y$-directions but differently in $t$-direction. TGV-$\LLambda$ has a similar U-Net architecture but with 4 output channels instead of 2. Consequently, for TGV-$\LLambda$, the estimated parameter maps take the form 
$\LLambda_{1, \Theta}:=[\LLambda_{1, \Theta}^{t}, \LLambda_{1, \Theta}^{xy}, \LLambda_{1, \Theta}^{xy}]^\trans$ and $\LLambda_{0, \Theta}:=[\LLambda_{0, \Theta}^{t}, \LLambda_{0, \Theta}^{xy}, \LLambda_{0, \Theta}^{xy}]^\trans$. 
The number of trainable parameters for the TV-$\LLambda$ and TGV-$\LLambda$ networks is  83,986 and 84,020, respectively.
The number of unrolled iterations for TV-$\LLambda$ and TGV-$\LLambda$ is, again, $T=96$ and $T=128$, respectively.

\paragraph{Training}
For supervised training, we randomly split 108 images from the MDCNN dataset \cite{mdcnn} into 80 images for training, 8 images for validation, and 20 images for testing. The same test set is used for zero-shot self-supervised learning for a direct comparison. 
Together with the U-Net parameters, we train an additional parameter \(t\). The output of the U-Net is passed through the softplus functions with $\beta=1$,
then multiplied with softplus(\(t\)) (also with $\beta=1$) to obtain the final regularization parameter maps. The parameter \(t\) is initialized to 0 before training. Adding this additional trainable parameter was observed to help the training converge more quickly. Each model was trained for 100 epochs. For both supervised and zero-shot training, we use the Adam optimizer with a learning rate of $2 \times 10^{-4}$. For the SSDU training, at each epoch we split the index set  $J:=I_{\mathrm{Rec}}\uplus I_{\mathrm{Loss}} \uplus I_{\mathrm{Val}}$ by randomly assigning the coefficients of entire $k$-space lines to the sets $I_{\mathrm{Rec}}, I_{\mathrm{Loss}}$ and $I_{\mathrm{Loss}}$. 
We use a split of $81\%/9\%/10$ and a batch size of 1.
Since training the TGV-$\LLambda$ model normally would overflow our 32GB GPU memory capacity,
we applied the gradient checkpoint technique \cite{chen2016memory}
to reduce memory usage at the cost of runtime.

\section{Results}

\subsection{Image Denoising Results}

Table \ref{tab:denoising_table} quantitatively summarizes the obtained results in terms of PSNR, MSE, SSIM, and HaarPSI by listing the means and standard deviations of the respective metrics calculated over the test dataset. As we can see, the general tendency is that the supervised methods perform best, as expected. This observation holds for RAM, TV-$\LLambda$ and TGV-$\LLambda$. 
Where our models differ is that the gap between the performance of the models trained in a supervised manner and the ones trained in a zero-shot fashion is much smaller for TV-$\LLambda$ and TGV-$\LLambda$ compared to RAM. 
Further, we see that RAM tends to suffer from overfitting for the ZS-TTT case, especially at higher noise levels, since there are notable differences between performances obtained by running the ZS-TTT of RAM training until the last iteration and stopping appropriately, interrupting the training procedure. 
For example, for $\sigma=0.15$, RAM trained by ZS-TTT reaches an SSIM statistic of $0.27 \pm 0.14 $ when run for 500 iterations, and an SSIM of $0.64 \pm 0.04$ when ZS-TTT is stopped according to the best validation SSDU loss. 
In contrast, for the same noise level, for TV-$\LLambda$ and TGV-$\LLambda$, there is almost no difference between ZS-TTT run for 500 iterations and ZS-TTT stopped according to the best validation SSDU loss, where they reach SSIM values of $0.70 \pm 0.05$ and $0.71 \pm 0.05$ (TV-$\LLambda$), and  $0.72 \pm 0.06$ and $0.72 \pm 0.06 $ (TGV-$\LLambda$), respectively.

\begin{table}[]
    \centering
    \fontsize{7.5pt}{8pt}\selectfont
        \renewcommand{\arraystretch}{1.15}  
\begin{tabular}{c c c c c c c}
                                                                &                                                & \bf{Strat.}                  & \bf{Stop}     & \bf{PSNR} $\uparrow$ & \bf{SSIM} $\uparrow$ & \bf{HaarPSI} $\uparrow$ \\
    \Thickline
    \multirow{21}{*}{\rotatebox[origin=c]{90}{$\sigma = 0.05$}} & \multirow{7}{*}{\rotatebox[origin=c]{90}{RAM}} & \textit{Sup.}                & N.A.          & $32.87 \pm 2.75$     & $0.91 \pm 0.03$      & $0.86 \pm 0.02$         \\  
    \cline{3-7}
                                                                &                                                & \multirow{3}{*}{\textit{ST}} & \textit{Last} & $27.55 \pm 2.57$     & $0.77 \pm 0.07$      & $0.70 \pm 0.06$         \\  
                                                                &                                                &                              & \textit{Val.} & $31.48 \pm 4.04$     & $0.91 \pm 0.03$      & $0.85 \pm 0.04$         \\  
                                                                &                                                &                              & \textit{Test} & $32.97 \pm 2.83$     & $0.92 \pm 0.03$      & $0.86 \pm 0.02$         \\
\end{tabular}

        \\[0.2em]
    \caption{
        \footnotesize
        Result summary for image denoising.
        \ ``\textbf{Strat.}'' is the strategy, where ``\textit{Sup.}'', ``\textit{ST}'', and ``\textit{ZS}'' 
        denote the supervised, supervised + TTT, and zero-shot TTT strategies, respectively.
        \ ``\textbf{Stop}'' is the stopping criterion,
        where ``\textit{Last}'', ``\textit{Val.}'', and ``\textit{Test}'' denote the last iteration, 
        the iteration with the lowest validation SSDU loss, 
        and the iteration with the lowest test MSE (oracle), respectively.
    }
    \label{tab:denoising_table}
\end{table}

\begin{table}[ht]
    \centering
    \begin{tabular}{l|c c c }
                             & RAM      & TV-$\LLambda$ & TGV-$\LLambda$ \\
        \hline
        Inference time [min] & $\sim 2$ & $\sim 2$      & $\sim 5$       \\
    \end{tabular}\\[0.5em]
    \caption{Zero-shot self-supervised inference times for RAM, TV-$\LLambda$ and TGV-$\LLambda$. The reported times correspond to 500 model weight updates on a GPU with 24GB memory.
    }
    \label{tab:denoising_times}
\end{table}

\begin{figure}[]
    \centering
    \includegraphics[width=\linewidth]{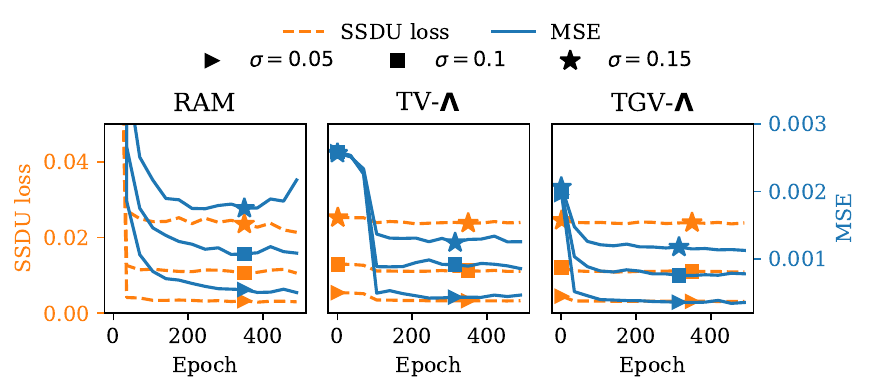}\\[-1em]
    \caption{Training behaviour of the self-supervised zero-shot test-time training (ZS-TTT) process for RAM, TV-$\LLambda$ and TGV-$\LLambda$ for the image denoising example. TV-$\LLambda$ and TGV-$\LLambda$ surpass RAM in terms of achievable MSE with respect to the target image. Further, while TV-$\LLambda$ and TGV-$\LLambda$ do not suffer from overfitting, overfitting is visible for RAM for the highest noise level.}
    \label{fig:denoising_ssdu_loss_vs_mse}
\end{figure}

Figure \ref{fig:denoising_ssdu_loss_vs_mse} shows the comparison of the SSDU loss used for ZS-TTT and the MSE with respect to the target image for an example of the image denoising problem, 
which illustrates the behaviour of the methods during the optimization process. First, we notice how, in general, the optimization of RAM seems to decrease the achievable MSE more slowly, while TV-$\LLambda$ and TGV-$\LLambda$ reach a rather good MSE value after $\sim$100 iterations of self-supervised training. Again, we see the absence of overfitting for both TV-$\LLambda$ and TGV-$\LLambda$, also for the highest considered noise level, where RAM, in contrast, overfits to the noise.

Table \ref{tab:denoising_times} compares the times required for zero-shot self-supervised learning of the models, indicating how quickly RAM, TV-$\LLambda$ and TGV-$\LLambda$ can learn to denoise an image from scratch.

Figure \ref{fig:denoising_recons} shows an example of a comparison of RAM, TV-$\LLambda$ and TGV-$\LLambda$ for the three considered setups, i.e., supervised training, supervised + TTT, and ZS-TTT for $\sigma=0.1$. The figure visually confirms what was previously observed in Table \ref{tab:denoising_table}, where we see that, as expected, the supervised models perform best. However, as previously noted, the performance gap between the supervised and the ZS-TTT TV-$\LLambda$ and TGV-$\LLambda$ is much smaller than for RAM. Figure \ref{fig:denoising_parameter_maps} compares the regularization parameter maps corresponding to the example in Figure \ref{fig:denoising_recons}.
Interestingly, we see that the regularization parameter maps are relatively different from each other. More precisely, for the supervised and supervised + TTT-finetuning case, the regularization parameter maps seem to contain the noise realization that is present in the noisy image, while the ones obtained by ZS-TTT are much smoother.

\newcommand{\denoiseroot}{images/denoise}  

\newlength{\denoisecellwidth}
\setlength{\denoisecellwidth}{0.15\textwidth}

\newlength{\denoiserowlabelheight}
\setlength{\denoiserowlabelheight}{0.9\denoisecellwidth} 

\newcommand{\denoiserotaterowlabel}[1]{%
  \rotatebox[origin=c]{90}{ \makebox[\denoiserowlabelheight][c]{#1} }
}

\newcommand{\imgrelcoord}[3]{%
  ($($(#1.south west)!#2!(#1.south east)$)!#3!($(#1.north west)!#2!(#1.north east)$)$)%
}

\newlength{\denoisespysize}

\NewDocumentCommand{\denoiseimagezoombase}{m m O{} O{} O{} O{0.2} O{0.28} O{0.22} O{0.78}}{%
  \setlength{\denoisespysize}{0.35\linewidth}%

  \begin{tikzpicture}[
      spy using outlines={shape=rectangle,
          magnification=3,
          size=\denoisespysize,
          draw=#2,
          connect spies
        }
    ]
    \node[inner sep=0] (img) {%
      \fbox{\includegraphics[width=0.9\linewidth]{#1}}%
    };

    \spy on \imgrelcoord{img}{#6}{#7}
      in node at \imgrelcoord{img}{#8}{#9};

    \ifstrempty{#3}{}{%
      \node[
        anchor=south west,
        font=\tiny\bfseries,
        text=white,
        fill=black,
        fill opacity=0.55,
        text opacity=1,
        inner sep=1pt
      ] at ([xshift=1pt,yshift=1pt]img.south west) {SSIM: #3};
    }%

    \ifstrempty{#4}{}{%
      \node[
        anchor=south east,
        font=\tiny\bfseries,
        text=white,
        fill=black,
        fill opacity=0.55,
        text opacity=1,
        inner sep=1pt
      ] at ([xshift=-1pt,yshift=1pt]img.south east) {PSNR: #4};
    }%

    \ifstrempty{#5}{}{%
      \node[
        anchor=north east,
        font=\tiny\bfseries,
        text=white,
        fill=black,
        fill opacity=0.55,
        text opacity=1,
        inner sep=1pt
      ] at ([xshift=-1pt,yshift=-1pt]img.north east) {HaarPSI: #5};
    }%
  \end{tikzpicture}%
}

\NewDocumentCommand{\denoiseimagezoom}{m O{} O{} O{} O{white}}{%
  \denoiseimagezoombase{\denoiseroot/#1.pdf}{#5}[#2][#3][#4]%
}

\NewDocumentCommand{\denoiseimagezoomsub}{m m O{} O{} O{} O{white}}{%
  \begin{minipage}[c]{#2}
    \centering
    \denoiseimagezoom{#1}[#3][#4][#5][#6]%
  \end{minipage}%
}

\NewDocumentCommand{\denoisereconsub}{m m O{} O{} O{} O{white}}{%
  \denoiseimagezoomsub{#1/reconstruction}{\denoisecellwidth}[#3][#4][#5][#6]
}

\NewDocumentCommand{\denoiserecons}{m m m O{white}}{
  \begin{figure}[]
    \centering
    \scriptsize
    \setlength{\tabcolsep}{1pt}
    \renewcommand{\arraystretch}{1.0}
    \captionsetup[subfigure]{font=scriptsize,justification=centering}

    \begin{tabular}{cccc}

       & Ground-truth
       & Noisy
       &
      \\[-0.12em]

       & \denoiseimagezoomsub{#1/target_image}{\denoisecellwidth}[][][][#4]
       & \denoiseimagezoomsub{#1/noisy_image}{\denoisecellwidth}[0.23][20.02][0.40][#4]
       &
      \\
      \noalign{\vskip 0.2em}

      \hline                                                        \\[-0.8em]

       & Supervised
       & Supervised + TTT
       & ZS-TTT
      \\[-0.18em]

      \denoiserotaterowlabel{RAM}
       & \denoisereconsub{#1/ram/sup}{}[0.89][32.30][0.76][#4]     
       & \denoisereconsub{#1/ram/supttt_best_val}{}[0.90][32.43][0.77][#4]  
       & \denoisereconsub{#1/ram/zs_best_val
         }{}[0.73][29.40][0.65][#4]             
      \\[-0.1em]

      \denoiserotaterowlabel{TV-$\LLambda$}
       & \denoisereconsub{#1/tv/sup}{}[0.89][32.05][0.76][#4]      
       & \denoisereconsub{#1/tv/supttt_best_val}{}[0.89][31.79][0.74][#4]  
       & \denoisereconsub{#1/tv/zs_final}{}[0.79][30.37][0.79][#4]          
      \\[-0.1em]

      \denoiserotaterowlabel{TGV-$\LLambda$}
       & \denoisereconsub{#1/tgv/sup}{}[0.89][32.11][0.76][#4]     
       & \denoisereconsub{#1/tgv/supttt_best_val}{}[0.89][31.99][0.76][#4] 
       & \denoisereconsub{#1/tgv/zs_final}{}[0.86][31.26][0.71][#4]         
      \\[-0.1em]
    \end{tabular}
    \caption{#3}
    \label{#2}
  \end{figure}
}

\denoiserecons{sample_552_sigma_0_1_magma}{fig:denoising_recons}{
  Reconstructions with different denoising methods on a test sample with noise standard deviation $\sigma = 0.1$.
  Metrics for ZS-TTT (zero-shot) TV-$\LLambda$ and TGV-$\LLambda$ are computed on the final reconstructions (after 500 epochs),
  while others are computed on the best validation-SSDU-loss epoch.
}[green]

\newlength{\denoisemapcellwidth}
\setlength{\denoisemapcellwidth}{0.96\denoisecellwidth}

\NewDocumentCommand{\denoisemapsub}{m m O{} O{} O{white}}{%
  \denoiseimagezoomsub{#1}{\denoisemapcellwidth}[][][][#3]%
}

\NewDocumentCommand{\denoisecolorbar}{m}{%
  \begin{minipage}[c]{0.3\denoisemapcellwidth}
    \centering
    \rotatebox[origin=c]{90}{
    \includegraphics[width=0.88\denoisemapcellwidth]{\denoiseroot/#1_colorbar.pdf}%
    }
  \end{minipage}%
}

\NewDocumentCommand{\denoisemaps}{m m m O{white}}{
  \begin{figure}[]
    \centering
    \scriptsize
    \setlength{\tabcolsep}{0pt}
    \renewcommand{\arraystretch}{1.0}
    \captionsetup[subfigure]{font=scriptsize,justification=centering}

    \begin{tabular}{ccccc}

       & Supervised
       & Supervised + TTT
       & ZS-TTT
       &
      \\[-0.15em]

      \denoiserotaterowlabel{TV $\LLambda$-map}
       & \denoisemapsub{#1/tv/sup/lambda_1_v}{}[#4]
       & \denoisemapsub{#1/tv/supttt_best_val/lambda_1_v}{}[#4]
       & \denoisemapsub{#1/tv/zs_final/lambda_1_v}{}[#4]
       & \denoisecolorbar{#1/tv/zs_final/lambda_1_v}
      \\[-0.16em]

      \denoiserotaterowlabel{TGV $\LLambda_1$-map}
       & \denoisemapsub{#1/tgv/sup/lambda_1_v}{}[#4]
       & \denoisemapsub{#1/tgv/supttt_best_val/lambda_1_v}{}[#4]
       & \denoisemapsub{#1/tgv/zs_final/lambda_1_v}{}[#4]
       & \denoisecolorbar{#1/tgv/zs_final/lambda_1_v}
      \\[-0.16em]

      \denoiserotaterowlabel{TGV $\LLambda_0$-map}
       & \denoisemapsub{#1/tgv/sup/lambda_0_w}{}[#4]
       & \denoisemapsub{#1/tgv/supttt_best_val/lambda_0_w}{}[#4]
       & \denoisemapsub{#1/tgv/zs_final/lambda_0_w}{}[#4]
       & \denoisecolorbar{#1/tgv/zs_final/lambda_0_w}
      \\[-0.16em]

      \denoiserotaterowlabel{TGV $\text{log}_{10}(\LLambda_0 / \LLambda_1)$}
       & \denoisemapsub{#1/tgv/sup/lambda_0_w_over_lambda_1_v_log10}{}[#4]
       & \denoisemapsub{#1/tgv/supttt_best_val/lambda_0_w_over_lambda_1_v_log10}{}[#4]
       & \denoisemapsub{#1/tgv/zs_final/lambda_0_w_over_lambda_1_v_log10}{}[#4]
       & \denoisecolorbar{#1/tgv/zs_final/lambda_0_w_over_lambda_1_v_log10}
      \\[-0.12em]
    \end{tabular}
    \caption{#3}
    \label{#2}
  \end{figure}
}

\denoisemaps{sample_552_sigma_0_1_magma}{fig:denoising_parameter_maps}{
  Regularization parameter maps produced by our TV-$\LLambda$ and TGV-$\LLambda$ methods for supervised training, supervised training + test-time training and zero-shot test-time training. While the parameter maps estimated by network trained in a supervised fashion, noise is present in the $\LLambda$-maps, the ones estimated by the ZS-TTT network are much smoother.  Note that if the ratio $\LLambda_0/\LLambda_1$ shown at the last row is high, one expects TGV to behave roughly like TV \cite{Papafitsoros_Valkonen_2015}.
}[green]

\subsection{Dynamic Cardiac MRI Results}

\begin{figure}[]
    \centering
    \includegraphics[width=\linewidth]{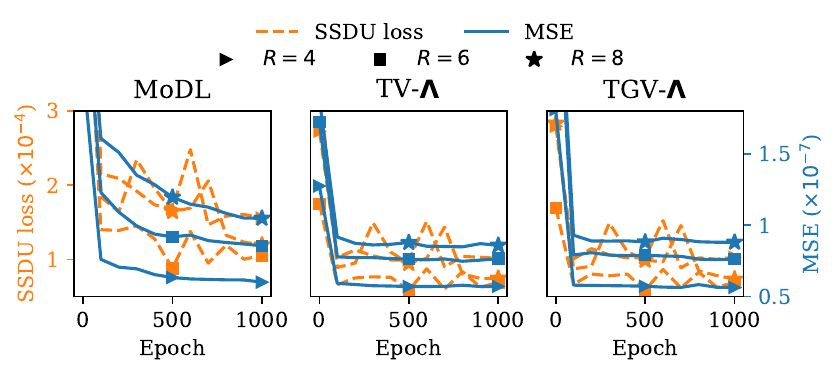}\\[-1em]
    \caption{Training behaviour of the self-supervised zero-shot test-time training (ZS-TTT) process for MoDL, TV-$\LLambda$ and TGV-$\LLambda$ for the dynamic cardiac MR imaging problem. MoDL, TV-$\LLambda$ and TGV-$\LLambda$ yield similar performance. However, TV-$\LLambda$ and TGV-$\LLambda$ nearly achieve this performance already after 100-200 iterations, while MoDL requires substantially more iterations to reach a similar MSE.}
    \label{fig:mri_ssdu_loss_vs_mse}
\end{figure}

Table \ref{tab:mri_table} quantitatively summarizes the obtained results in terms of PSNR, MSE, SSIM, and HaarPSI by listing the mean and standard deviations of the respective metrics calculated over the test dataset. Again, the supervised methods tend to perform best, as expected. In contrast to the denoising experiment, we see that for TV-$\LLambda$ and TGV-$\LLambda$, the difference between the supervised models and the ones obtained by ZS-TTT is often negligible. For example, even for the largest acceleration factor $R=6$, TV-$\LLambda$ yields SSIM values of $0.83 \pm 0.04$ (supervised) and $0.82 \pm 0.04$ (ZS-TTT; last iteration), and TGV-$\LLambda$ also yields similar SSIM values of $0.83 \pm 0.04$ (supervised) and $0.82 \pm 0.04$ (ZS-TTT; last iteration). In contrast, MoDL yields, for example, SSIM values $0.82 \pm 0.04$ (supervised) and $0.73 \pm 0.05$ (ZS-TT; best validation SSDU loss). This observation holds across the different image metrics, with MoDL consistently exhibiting a larger gap in performance between the supervised setting and the ZS-setting, while TV-$\LLambda$ and TGV-$\LLambda$ show only relatively small differences.

\begin{table}[]
    \centering
    \fontsize{7.5pt}{8pt}\selectfont
        \renewcommand{\arraystretch}{1.15}  
\begin{tabular}{c c c c c c c}
                                                            &                                                 & \bf{Strat.}                  & \bf{Stop}     & \bf{PSNR} $\uparrow$ & \bf{SSIM} $\uparrow$ & \bf{HaarPSI} $\uparrow$ \\
        \Thickline
        \multirow{21}{*}{\rotatebox[origin=c]{90}{$R = 4$}} & \multirow{7}{*}{\rotatebox[origin=c]{90}{MoDL}} & \textit{Sup.}                & N.A.          & $36.09 \pm 1.74$     & $0.86 \pm 0.04$      & $0.89 \pm 0.03$         \\  
        \cline{3-7}
                                                            &                                                 & \multirow{3}{*}{\textit{ST}} & \textit{Last} & $35.94 \pm 1.67$     & $0.85 \pm 0.04$      & $0.87 \pm 0.03$         \\  
                                                            &                                                 &                              & \textit{Val.} & $36.09 \pm 1.74$     & $0.86 \pm 0.04$      & $0.89 \pm 0.03$         \\  
                                                            &                                                 &                              & \textit{Test} & $36.10 \pm 1.74$     & $0.86 \pm 0.04$      & $0.88 \pm 0.03$         \\
\end{tabular}

        \\[0.2em]
    \caption{
        \footnotesize
        Result summary for the dynamic MR image reconstruction problem.
        \ ``\textbf{Strat.}'' is the strategy, where ``\textit{Sup.}'', ``\textit{ST}'', and ``\textit{ZS}'' 
        denote the supervised, supervised + TTT, and zero-shot TTT strategies, respectively.
        \ ``\textbf{Stop}'' is the stopping criterion,
        where ``\textit{Last}'', ``\textit{Val.}'', and ``\textit{Test}'' denote the last iteration, 
        the iteration with the lowest validation SSDU loss, 
        and the iteration with the lowest test MSE (oracle), respectively.
    }
    \label{tab:mri_table}
\end{table}

\begin{table}[ht]
    \centering
    \begin{tabular}{l|c c c }
                             & MoDL     & TV-$\LLambda$ & TGV-$\LLambda$ \\
        \hline
        Inference time [min] & $\sim 4$ & $\sim 2$      & $\sim 7$       \\
    \end{tabular}\\[0.5em]
    \caption{Zero-shot self-supervised inference times for MoDL, TV-$\LLambda$ and TGV-$\LLambda$. 
     The reported times for ZS-TTT correspond to 100 model weight updates on a GPU with 32GB memory.
    }
    \label{tab:mri_times}
\end{table}

\newcommand{\mriroot}{images/mri}  

\newlength{\mricellwidth}
\setlength{\mricellwidth}{0.11\textwidth}

\newlength{\mrirowlabelheight}
\setlength{\mrirowlabelheight}{\mricellwidth} 

\newcommand{\mrirotaterowlabel}[1]{%
  \rotatebox[origin=c]{90}{ \makebox[\mrirowlabelheight][c]{#1} }
}

\newcommand{\mridashedline}[2]{
  \draw[
    purple,
    line width=0.1em, draw opacity=0.6,
    dashed
  ] #1 -- #2;
}

\newcommand{\mrisolidline}[2]{
  \draw[
    blue,
    line width=0.1em, draw opacity=0.6,
      -{Stealth[length=5pt,width=2pt]}
  ] #1 -- #2;
}

\newcommand{\mrimetricfont}{\fontsize{4.1pt}{4.11pt}\bfseries}

\NewDocumentCommand{\mriimagezoombase}{m m O{} O{} O{} O{\mricellwidth}}{%
  \begin{tikzpicture}[
    ]

    \node[inner sep=0] (img) {\includegraphics[width=#6]{#1}};


    \mrisolidline{([xshift=0.3pt]img.west)}{([xshift=-0.3pt]img.east)}

    \mridashedline{([yshift=-0.3pt]img.north)}{([yshift=0.3pt]img.south)}

    \ifstrempty{#3}{}{%
      \node[
        anchor=south west,
        font=\mrimetricfont,
        text=white,
        fill=black,
        fill opacity=0.55,
        text opacity=1,
        inner sep=1pt
      ] at ([xshift=1pt,yshift=1pt]img.south west) {SSIM: #3};
    }%

    \ifstrempty{#4}{}{%
      \node[
        anchor=south east,
        font=\mrimetricfont,
        text=white,
        fill=black,
        fill opacity=0.55,
        text opacity=1,
        inner sep=1pt
      ] at ([xshift=-1pt,yshift=1pt]img.south east) {PSNR: #4};
    }%

    \ifstrempty{#5}{}{%
      \node[
        anchor=north east,
        font=\mrimetricfont,
        text=white,
        fill=black,
        fill opacity=0.55,
        text opacity=1,
        inner sep=1pt
      ] at ([xshift=-1pt,yshift=-1pt]img.north east) {HaarPSI: #5};
    }%

  \end{tikzpicture}%
}

\NewDocumentCommand{\mriimagezoomsubmetrics}{m m O{} O{} O{} O{\mricellwidth}}{%
  \ifstrempty{#2}{%
    \begin{minipage}[t]{#6}
      \centering
      \mriimagezoombase{\mriroot/#1.pdf}{white}[#3][#4][#5][#6]
    \end{minipage}%
  }{%
    \subfloat[#2\label{fig:mri:#1}]{%
      \begin{minipage}[t]{#6}
        \centering
        \mriimagezoombase{\mriroot/#1.pdf}{white}[#3][#4][#5][#6]
      \end{minipage}%
    }%
  }%
}

\NewDocumentCommand{\viewxt}{m O{\mricellwidth}}{
  \begin{tikzpicture}
    \node[inner sep=0] (img) {\includegraphics[height=#2]{\mriroot/#1.pdf}};
    \mridashedline{([yshift=-0.3pt, xshift=0.3pt]img.north west)}{([yshift=0.3pt, xshift=0.3pt]img.south west)}
  \end{tikzpicture}
}

\NewDocumentCommand{\viewyt}{m O{\mricellwidth}}{
  \begin{tikzpicture}
    \node[inner sep=0] (img) {\includegraphics[width=#2, angle=90]{\mriroot/#1.pdf}};
    \mrisolidline{([yshift=0.3pt, xshift=0.3pt]img.south west)}{([yshift=-0.3pt, xshift=0.3pt]img.north west)}
  \end{tikzpicture}
}

\NewDocumentCommand{\mripanel}{m m m O{} O{} O{} O{\mricellwidth}}{%
  \viewxt{#1_view_x_time_col_103_#2}[#7]\hspace{-1.0em}
  \viewyt{#1_view_y_time_row_103_#2}[#7]\hspace{-0.6em}
  \mriimagezoomsubmetrics{#1_frame_00_#2}{#3}[#4][#5][#6][#7]%
}

\NewDocumentCommand{\mripanelitem}{m m m m O{} O{} O{} O{\mricellwidth}}{%
  \begin{tabular}{@{}c@{}}%
    \mripanel{#1}{#2}{#3}[#5][#6][#7][#8] \\[#4]%
  \end{tabular}%
}

\NewDocumentCommand{\mriimagezoomsub}{m m O{} O{} O{} O{\mricellwidth}}{%
  \mripanelitem{#1}{vmin_0.0e+00_vmax_6.0e-03}{#2}{0em}[#3][#4][#5][#6]%
}


\readlist\metricssense{0.57,29.01,0.55} 
\readlist\metricsmodlsup{0.81,34.23,0.79}
\readlist\metricsmodlsupttt{0.79,33.59,0.80}  
\readlist\metricsmodlzs{0.74,32.50,0.74}  
\readlist\metricstvsup{0.80,33.97,0.80}
\readlist\metricstvsupttt{0.80,34.07,0.84}  
\readlist\metricstvzs{0.79,33.79,0.82}  
\readlist\metricstgvsup{0.80,34.17,0.82}
\readlist\metricstgvsupttt{0.79,34.03,0.84}  
\readlist\metricstgvzs{0.79,33.75,0.82}  

\NewDocumentCommand{\mrirecons}{m m}{
  \begin{figure}[]
    \centering
    \scriptsize
    \setlength{\tabcolsep}{1pt}
    \renewcommand{\arraystretch}{1.0}
    \captionsetup[subfigure]{font=scriptsize,justification=centering}

    \begin{tabular}{cccc}

       & Ground-truth
       & Iterative SENSE
       &
      \\[-0.12em]

       & \mriimagezoomsub{#1/ground_truth}{}[][]
       &
      \mriimagezoomsub{#1/sense/reconstruction}{}[\metricssense[1]][\metricssense[2]][\metricssense[3]]
       &
      \\[-0.1em]

      \hline                                                                  \\[-0.8em]

       & Supervised
       & Supervised + TTT
       & ZS-TTT
      \\[-0.12em]

      \mrirotaterowlabel{MoDL}
       & \mriimagezoomsub{#1/modl/sup/reconstruction}{}[\metricsmodlsup[1]][\metricsmodlsup[2]][\metricsmodlsup[3]]
       & \mriimagezoomsub{#1/modl/supttt_best_val/reconstruction}{}[\metricsmodlsupttt[1]][\metricsmodlsupttt[2]][\metricsmodlsupttt[3]]
       & \mriimagezoomsub{#1/modl/zs_final/reconstruction}{}[\metricsmodlzs[1]][\metricsmodlzs[2]][\metricsmodlzs[3]]
      \\[-0.1em]

      \mrirotaterowlabel{TV-$\LLambda$}
       & \mriimagezoomsub{#1/tv/sup/reconstruction}{}[\metricstvsup[1]][\metricstvsup[2]][\metricstvsup[3]]           
       & \mriimagezoomsub{#1/tv/supttt_best_val/reconstruction}{}[\metricstvsupttt[1]][\metricstvsupttt[2]][\metricstvsupttt[3]]
       & \mriimagezoomsub{#1/tv/zs_final/reconstruction}{}[\metricstvzs[1]][\metricstvzs[2]][\metricstvzs[3]]
      \\[-0.1em]

      \mrirotaterowlabel{TGV-$\LLambda$}
       & \mriimagezoomsub{#1/tgv/sup/reconstruction}{}[\metricstgvsup[1]][\metricstgvsup[2]][\metricstgvsup[3]]             
       & \mriimagezoomsub{#1/tgv/supttt_best_val/reconstruction}{}[\metricstgvsupttt[1]][\metricstgvsupttt[2]][\metricstgvsupttt[3]]
       & \mriimagezoomsub{#1/tgv/zs_final/reconstruction}{}[\metricstgvzs[1]][\metricstgvzs[2]][\metricstgvzs[3]]
      \\[-0.1em]
    \end{tabular}
    \caption{#2}
    \label{fig:mri_recons}
  \end{figure}
}

\mrirecons{sample_5_R6_magma}{
  Reconstructions with different methods on a test sample with acceleration 
  $R = 6$.
  Metrics for ZS-TTT (zero-shot) are computed on the final reconstructions (after 100 epochs),
  while others are computed on the best validation-SSDU-loss epoch.
}

\newlength{\mrilambdacellwidth}
\setlength{\mrilambdacellwidth}{0.1\textwidth}

\readlist\tvxy{lambda_map_idx1_spatial,0.0e+00,1.0e-01}
\readlist\tvt{lambda_map_idx0_temporal,0.0e+00,4.0e-01}
\readlist\tgvxyone{lambda_map_idx1_spatial_1,0.0e+00,1.0e-01}
\readlist\tgvxyzero{lambda_map_idx4_spatial_0,0.0e+00,2.0e-01}
\readlist\tgvxyratio{lambda_0_w_over_lambda_1_v_spatial_log10,-2.0e+00,2.0e+00}
\readlist\tgvtone{lambda_map_idx0_temporal_1,0.0e+00,1.0e+00}
\readlist\tgvtzero{lambda_map_idx3_temporal_0,0.0e+00,2.0e-01}
\readlist\tgvtratio{lambda_0_w_over_lambda_1_v_temporal_log10,-2.0e+00,2.0e+00}

\NewDocumentCommand{\mrilambdacell}{m m m}{%
  \mripanelitem{#1/lambda_maps/#2[1]}{vmin_#2[2]_vmax_#2[3]}{#3}{0em}[][][][\mrilambdacellwidth]%
}

\NewDocumentCommand{\mricolorbar}{m m}{%
  \begin{minipage}[c]{0.3\mrilambdacellwidth}
    \centering
    \rotatebox[origin=c]{90}{
      \includegraphics[width=\mrilambdacellwidth]{\mriroot/#1/lambda_maps/#2[1]_evolution_vmin_#2[2]_vmax_#2[3]_colorbar.pdf}%
    }
  \end{minipage}%
}

\NewDocumentCommand{\mrimaps}{m m}{
  \begin{figure}[]
    \centering
    \scriptsize
    \setlength{\tabcolsep}{1pt}
    \renewcommand{\arraystretch}{1.0}
    \captionsetup[subfigure]{font=scriptsize,justification=centering}

    \begin{tabular}{ccccc}

       & Supervised
       & Supervised + TTT
       & ZS-TTT
       &
      \\[-0.1em]

      \mrirotaterowlabel{TV $\LLambda^{xy}$-map}
       & \mrilambdacell{#1/tv/sup}{\tvxy}{}
       & \mrilambdacell{#1/tv/supttt_best_val}{\tvxy}{}
       & \mrilambdacell{#1/tv/zs_final}{\tvxy}{}
       & \mricolorbar{#1/tv/zs_final}{\tvxy}
      \\[-0.1em]

      \mrirotaterowlabel{TV $\LLambda^{t}$-map}
       & \mrilambdacell{#1/tv/sup}{\tvt}{}
       & \mrilambdacell{#1/tv/supttt_best_val}{\tvt}{}
       & \mrilambdacell{#1/tv/zs_final}{\tvt}{}
       & \mricolorbar{#1/tv/zs_final}{\tvt}
      \\[-0.1em]

      \mrirotaterowlabel{TGV $\LLambda^{xy}_1$-map}
       & \mrilambdacell{#1/tgv/sup}{\tgvxyone}{}
       & \mrilambdacell{#1/tgv/supttt_best_val}{\tgvxyone}{}
       & \mrilambdacell{#1/tgv/zs_final}{\tgvxyone}{}
       & \mricolorbar{#1/tgv/zs_final}{\tgvxyone}
      \\[-0.1em]

      \mrirotaterowlabel{TGV $\LLambda^{t}_1$-map}
       & \mrilambdacell{#1/tgv/sup}{\tgvtone}{}
       & \mrilambdacell{#1/tgv/supttt_best_val}{\tgvtone}{}
       & \mrilambdacell{#1/tgv/zs_final}{\tgvtone}{}
       & \mricolorbar{#1/tgv/zs_final}{\tgvtone}
      \\[-0.1em]

      \mrirotaterowlabel{TGV $\LLambda^{xy}_0$-map}
       & \mrilambdacell{#1/tgv/sup}{\tgvxyzero}{}
       & \mrilambdacell{#1/tgv/supttt_best_val}{\tgvxyzero}{}
       & \mrilambdacell{#1/tgv/zs_final}{\tgvxyzero}{}
       & \mricolorbar{#1/tgv/zs_final}{\tgvxyzero}
      \\[-0.1em]

      \mrirotaterowlabel{TGV $\LLambda^{t}_0$-map}
       & \mrilambdacell{#1/tgv/sup}{\tgvtzero}{}
       & \mrilambdacell{#1/tgv/supttt_best_val}{\tgvtzero}{}
       & \mrilambdacell{#1/tgv/zs_final}{\tgvtzero}{}
       & \mricolorbar{#1/tgv/zs_final}{\tgvtzero}
      \\[-0.1em]

      \mrirotaterowlabel{TGV $\text{log}_{10}\left(\frac{\LLambda^{xy}_0}{\LLambda^{xy}_1}\right)$}
       & \mrilambdacell{#1/tgv/sup}{\tgvxyratio}{}
       & \mrilambdacell{#1/tgv/supttt_best_val}{\tgvxyratio}{}
       & \mrilambdacell{#1/tgv/zs_final}{\tgvxyratio}{}
       & \mricolorbar{#1/tgv/zs_final}{\tgvxyratio}
      \\[-0.1em]

      \mrirotaterowlabel{TGV $\text{log}_{10}\left(\frac{\LLambda^{t}_0}{\LLambda^{t}_1}\right)$}
       & \mrilambdacell{#1/tgv/sup}{\tgvtratio}{}
       & \mrilambdacell{#1/tgv/supttt_best_val}{\tgvtratio}{}
       & \mrilambdacell{#1/tgv/zs_final}{\tgvtratio}{}
       & \mricolorbar{#1/tgv/zs_final}{\tgvtratio}
    \end{tabular}
    \caption{#2}
    \label{fig:mri_parameter_maps}
  \end{figure}
}

\mrimaps{sample_5_R6_magma}{
  Regularization parameter maps produced by our TV-$\LLambda$ and TGV-$\LLambda$ methods for supervised training, supervised training + test-time training (Supervised + TTT), and zero-shot test-time training (ZS-TTT).
}

Figure \ref{fig:mri_recons} shows a comparison of the reconstructions obtained with all methods for all training settings. The figure visually aligns well with the observations visible from Table \ref{tab:mri_table}. The images obtained with ZS-TTT for TV-$\LLambda$ and TGV-$\LLambda$ are effectively indistinguishable from the ones obtained by supervised training. For MoDL, as can also be quantitatively seen in Table \ref{tab:mri_table}, this is not the case.

Figure \ref{fig:mri_parameter_maps} compares the regularization parameter maps obtained for TV-$\LLambda$ and TGV-$\LLambda$ across the three learning scenarios. Similarly to the denoising example, the regularization parameter maps are relatively different from each other, although the differences are less pronounced than for the image denoising case, especially for TV-$\LLambda$.

\section{Discussion}

In the following, we discuss the results obtained by DREAM as well as different aspects of the proposed combination of ZS-SSDU for learning the regularization maps for both TV-$\LLambda$ and TGV-$\LLambda$.

\subsection{Reparametrization of the Regularization Parameter Maps}

In \cite{kofler2023learning} and \cite{vu2025deep}, supervised learning approaches were considered to estimate deeply parametrized parameter maps avoiding the need to resort to time-consuming procedures such as bilevel optimization schemes, \cite{crockett2022bilevel}. One could indeed frame 
DREAM 
through the lenses of bilevel optimization by formulating the problem as:
\begin{equation}\label{eq:bilievel_problem}
    \begin{cases}
         & \underset{\LLambda}{\min} \sum_{(\YY, \XX_{\mathrm{true}}) \in \mathcal{D}} \mathcal{L}\big(\XX^{\ast}(\LLambda)\big), \\
         & \text{s.t.}\quad  \XX^{\ast}(\LLambda)\in \underset{\XX}{\arg\min} \, F_{\LLambda}(\XX),
    \end{cases}
\end{equation}
with $F_{\LLambda}$ denoting the functional of the lower-level problem, such as the adaptive TV or TGV functionals in \eqref{eq:tv_problem_Lambda} or \eqref{eq:tgv_problem_Lambda}, while $\mathcal{L}$ is an upper-level loss function. 
DREAM 
shares similarity with the formulation \eqref{eq:bilievel_problem} by simply setting $\mathcal{L}$ in \eqref{eq:bilievel_problem} as $\mathcal{L}_{\mathrm{ZS}}^{\mathrm{SSDU}}$ in \eqref{eq:loss_function_ssdu} and $\XX^{\ast}(\LLambda)$ as  the image $\XX_{\LLambda}^{T}$ obtained by the unrolled reconstruction scheme. 
However, one important difference is that in \eqref{eq:bilievel_problem}, one aims at directly estimating the adaptive regularization parameter maps $\LLambda$, while, as visible in \eqref{eq:loss_function_ssdu}, the loss function depends on $\Theta$ rather than $\LLambda$ directly.  The reason is that, in our case, the regularization parameter maps are parametrized as the output of the network $\mathrm{NET}_{\Theta}$. Such a reparametrization can, for example, also be found in another context, namely the deep image prior  framework \cite{ulyanov2018deep}, where, instead of optimizing over the image itself, the image is reparametrized as the output of a network applied to pure noise. There, the regularization mechanism is given by early stopping, the intuition being that both the neural network architecture and the optimizer chosen provide regularization for the sought image. On the other hand, in our case it is the reparametrization of the parameter maps $\LLambda_{\Theta}=\mathrm{NET}_{\Theta}(\XX_0)$ combined with the employed SSDU-loss that leads to meaningful and effective adaptive regularization parameter maps.

\begin{figure}[h!]
    \centering
    \includegraphics[width=\linewidth]{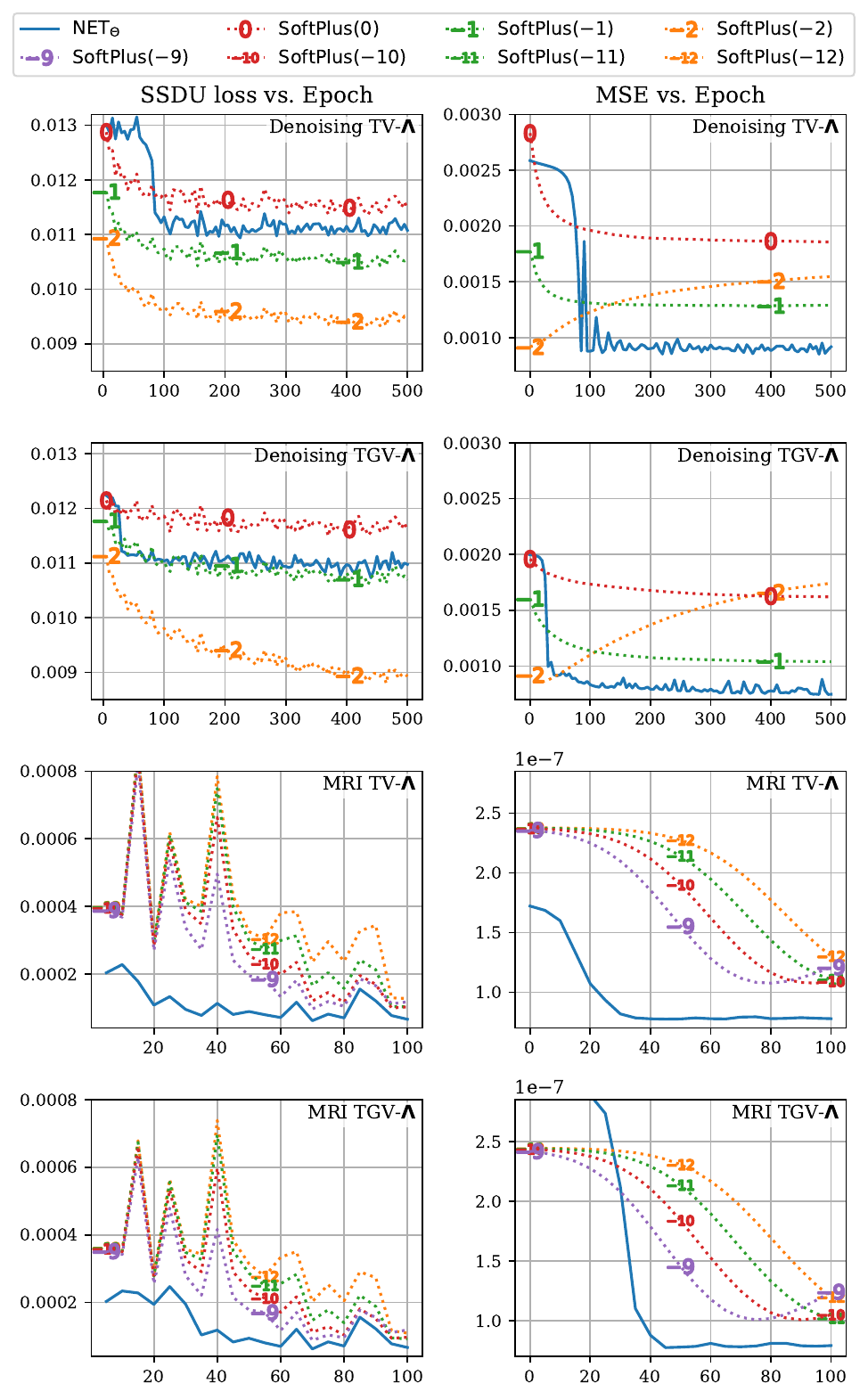}\\[-1em]
    \caption{
        ZS-SSDU loss and MSE for non-parametrized regularization-parameter maps in the denoising and MRI problems.
        Dotted lines indicate the explicit-map parametrization with different initial values. The reparametrization of the $\LLambda$-maps as output of the neural network, i.e., $\LLambda:=\mathrm{NET}_{\Theta}(\XX_0)$, allows to reach a much smaller MSE with respect to the target image.
    }
    \label{fig:explicit_map_ssdu_and_mse}
\end{figure}

Figure~\ref{fig:explicit_map_ssdu_and_mse} compares the ZS-SSDU loss and the MSE between the reconstruction and the target under two optimization strategies for TV-$\LLambda$ and TGV-$\LLambda$, for both image denoising and MRI reconstruction. In the first strategy, the entries of the regularization maps $\LLambda$ are treated directly as optimization variables. In the second, the maps are represented as
$
    \LLambda_{\Theta}
    =
    \mathrm{NET}_{\Theta}(\XX_{\mathrm{input}}),
$
and the optimization is performed over the network parameters $\Theta$. The latter strategy consistently yields substantially lower reconstruction MSE. In some cases, direct optimization over $\LLambda$ attains an even lower SSDU loss while producing a markedly higher target MSE, suggesting that the resulting maps overfit the self-supervised objective. This distinction is reminiscent of the deep image prior  framework. Here, however, the network represents the regularization maps rather than the reconstructed image. The resulting restriction of $\LLambda$ to maps that can be generated by $\mathrm{NET}_{\Theta}$ acts as an implicit prior, promoting structured parameter maps that generalize better than those obtained by unconstrained, direct optimization over their entries. This representation also provides an input-dependent reduction in the dimension of the optimization problem. In the denoising experiments, the network contains only $109{,}073$ parameters for TV-$\LLambda$ and $109{,}090$ for TGV-$\LLambda$, whereas direct optimization involves $512 \times 512 = 262{,}144$ and $512 \times 512 \times 2 = 524{,}288$ map entries, respectively. In the dynamic MRI experiments, the corresponding networks contain $83{,}986$ and $84{,}020$ parameters, compared with $208 \times 208 \times 25 \times 2 = 2{,}163{,}200$ and $208 \times 208 \times 25 \times 4 = 4{,}326{,}400$ entries in the directly optimized TV-$\LLambda$ and TGV-$\LLambda$ maps, respectively.

\subsection{Interpretability of the Regularization Maps}
We observe some structural differences between the regularization parameter maps inferred by the supervised and zero-shot approaches. 
For example, as seen from Figure \ref{fig:denoising_parameter_maps} in the denoising examples, 
the $\LLambda$ maps not only adapt to the image structure but also to the noise itself. 
This is not surprising since during training the model is provided with pairs of clean and noisy data and thus implicitly fed many instances of Gaussian noise,  eventually learning to recognize it and distinguish it from image structure. 
On the other hand, even though they still inherit the image structure, 
the zero-shot maps are flatter and more regular since the SSDU loss provides no noise information. 
Similar differences are also observed for the maps in the MRI experiment seen in Figure \ref{fig:mri_parameter_maps}, with the supervised ones exhibiting more details than the zero-shot ones, with the latter being smoother, however, these differences are not so pronounced as in denoising.
We hypothesize that the reason for the differences visible in the regularization parameter maps obtained by supervised and zero-shot self-supervised training essentially originates from the structure of the loss function in \eqref{eq:loss_function_ssdu}. In the recent review \cite{learning2026papafitsoros}, it is demonstrated that for image denoising, the network component $\mathrm{NET}_{\Theta}$ of TV-$\LLambda$ can learn to distinguish between image noise and image content, yielding parameter maps that are adapted not only to the noise level but also to the specific noise instance. Since no target images are used during the minimization of the SSDU loss function in \eqref{eq:loss_function_ssdu}, the network $\mathrm{NET}_{\Theta}$ is hindered from learning to discriminate between image content. Because the loss function \eqref{eq:loss_function_ssdu} involves the summation over different instances of measurement splits, its minimization yields an estimate of parameters $\Theta$ that is suitable ``on average`` and therefore yields smoother parameter maps. As also shown in \cite{learning2026papafitsoros}, for accelerated 2D MRI for TV-$\LLambda$, the learned parameter maps are much less adapted to the specific instances of noise and undersampling pattern compared to the image denoising setup. We believe this to be the reason that, in terms of features, the regularization parameter maps for the MRI example differ less from each other between the supervised and the zero-shot SSDU setting. 
Furthermore, for TV-$\LLambda$, the spatial and temporal components $\LLambda^{xy}$ and $\LLambda^{t}$ 
have complementary roles
in terms of their regularization contribution. In general, for the supervised, supervised + TTT, and the ZS-TTT setting, the achieved temporal regularization is much stronger than the spatial one,  except for regions depicting the cardiac movement, which is consistent with the results observed in \cite{kofler2023learning}. For TGV-$\LLambda$, on the other hand, the interplay between the spatial and temporal components is more complex due to the involvement of two parameters instead of one. For instance, as expected following the theory about the scalar versions \cite{Papafitsoros_Valkonen_2015}, it suffices that one of the parameters $\LLambda_{0}^{t}$, $\LLambda_{1}^{t}$ is locally small in order for the  temporal regularization effect of TGV to be locally small. From Figure \ref{fig:mri_parameter_maps} we observe indeed that it is $\LLambda_{0}^{t}$ which is locally small in the heart region for all three settings. However, the maps $\LLambda_{0}^{xy}$ and $\LLambda_{1}^{t}$ change significantly from the supervised to the ZS-TTT setting, which means that for TGV-$\LLambda$, dramatically different parameter combinations can lead to very similar reconstructions. We attribute that to the highly non-convex nature of the training problem.

\subsection{Performance, CPU Times and No Overfitting} 
From Figure \ref{fig:explicit_map_ssdu_and_mse}, we have seen how the achievable results also depend on the used reparametrization of the $\LLambda$-maps. Further, from Figures \ref{fig:denoising_ssdu_loss_vs_mse} and \ref{fig:mri_ssdu_loss_vs_mse} we have observed that for both considered problems, the deep learning-based models RAM and MoDL tend to exhibit a slower convergence to a set of suitable network parameters when trained in a zero-shot fashion using SSDU. For denoising, after 500 iterations, TV-$\LLambda$ and TGV-$\LLambda$ still yield a lower MSE with respect to the target image. For the dynamic MRI example, this effect is even more strongly visible, where for all acceleration factors $R=4,6,8$,  TV-$\LLambda$ and TGV-$\LLambda$ yield a lower MSE with respect to the target image after around 100 SSDU iterations than MoDL does after 1000 iterations, effectively yielding better reconstructions within  $\sim2$ minutes (TV-$\LLambda$) and $\sim7$ minutes (TGV-$\LLambda$) compared to $\sim4$ minutes (MoDL), see also Tables \ref{tab:denoising_table} and \ref{tab:mri_table}. 
Additionally, as can be seen from Table \ref{tab:mri_table}, for TV-$\LLambda$ and TGV-$\LLambda$, remarkably, the differences between the supervised models and the zero-shot models are very small, and sometimes even negligible, 
see also Figure \ref{fig:mri_recons}. 
In contrast, this is not the case for MoDL, which, when trained in a zero-shot fashion, substantially suffers in terms of achievable performance. While perhaps surprising at first glance, these results are plausible since MoDL, being a method for which the regularization is entirely learned from data, seems to profit more from being exposed to a larger dataset than being adapted to the specific unseen example. In contrast, for TV-$\LLambda$ and TGV-$\LLambda$, for which the regularization is hand-crafted, sample-specific adaptation of the regularization strength suffices to yield comparable performance to the supervised models. In addition to the observation that TV-$\LLambda$ and TGV-$\LLambda$ yield competitive results when compared to MoDL, this offers the possibility to rely on a simpler reconstruction method with improved interpretability without sacrificing image quality.
This important result is only somewhat analogously visible for the image denoising experiment. There, although the difference in performance between the supervised models and the zero-shot models is also smaller for TV-$\LLambda$ and TGV-$\LLambda$ compared to RAM, it is not possible to observe an effectively comparable performance of TV-$\LLambda$ and TGV-$\LLambda$. We believe that this can be attributed to the different nature of the two considered problems. For image denoising, there is inherently less data to be learned from for the ZS-TTT, which, in addition, is corrupted by potentially large noise. In contrast, for the MR example, there is more redundant data to be learned from due to the multi-coil setup of the problem, as can be seen from the forward model in \eqref{eq:multi_coil_mri_fwd}.

\subsection{Limitations}

Despite the advantages in terms of robustness with respect to overfitting and the competitive performance of the methods, 
DREAM has some limitations. 
First, the employed ZS-TTT framework for TV-$\LLambda$ and TGV-$\LLambda$ appears to suffer from longer reconstruction times. For example, as can be seen from Tables \ref{tab:denoising_times} and \ref{tab:mri_times}, inference can take up to two and five minutes for TV-$\LLambda$ and TGV-$\LLambda$, respectively, for a single 2D image, and up to two and seven minutes for TV-$\LLambda$ and TGV-$\LLambda$ for a single dynamic cardiac MR image. The reason for this, of course, mainly lies in the need to unroll the respective PDHG algorithms using a relatively large number of iterations $T$. Although lowering $T$ is possible, in \cite{kofler2023learning} it was noted that lowering $T$ during supervised training of TV-$\LLambda$ tends to yield too large regularization parameter maps to compensate for the smaller $T$, effectively leading to overly smoothed images. A promising direction that we aim to investigate in the future is to study different neural network architectures and initialization mechanisms to further accelerate the learning phase. 
We expect this strategy to be of interest for larger realistic reconstruction problems (e.g.\ 3D + time) or for other similar reconstruction methods such as the one considered in \cite{kofler2025ell1}, for which the considered problem is inherently larger due to its formulation using a synthesis-based approach.    

\section{Conclusions}
In this work, we presented a method for learning adaptive regularization parameter maps for both TV and TGV  regularization via algorithm unrolling in a zero-shot self-supervised framework. Through several examples, we have investigated the applicability of the approach for a simple 2D image denoising problem and for a challenging large-scale dynamic MR imaging problem (2D + time). 
In contrast to state-of-the-art models used for comparison (RAM for image denoising, MoDL for MRI), DREAM is fast, 
and it requires only a fraction of the weight updates 
as
the
methods of comparison to achieve a comparable performance. Furthermore, we demonstrated that representing the parameter maps as the output of a CNN applied to the input image provides a meaningful implicit prior and an input-dependent dimensionality reduction, substantially improving upon direct optimization of the regularization parameters. 
We have observed that the performance difference between supervised training and zero-shot test-time training is relatively small for TV-$\LLambda$ and TGV-$\LLambda$ for image denoising, and negligible for the dynamic cardiac MRI problem, while for RAM (for image denoising) and MoDL (for MRI), the differences are much larger. 
This implies that the learned TV-$\LLambda$ and TGV-$\LLambda$ regularization models can be used as valid alternatives to MoDL, with the advantage of exhibiting improved interpretability and, at the same time, not sacrificing image quality. 
Future work directions could include the extension of 
DREAM 
to other regularization functionals based on learned sparsifying transforms.

\bibliography{references}%

@InProceedings{Papafitsoros_Valkonen_2015,
author={Papafitsoros, K.
and Valkonen, T.},
editor={Aujol, Jean-Fran{\c{c}}ois
and Nikolova, Mila
and Papadakis, Nicolas},
title={Asymptotic Behaviour of Total Generalised Variation},
booktitle={Scale Space and Variational Methods in Computer Vision},
year={2015},
pages={702--714}
}

@article{Zou2023,
         author={Zou, Zihao
         and Liu, Jiaming
         and Wohlberg, Brendt 
         and Kamilov, Ulugbek S.},
         journal={IEEE Open J. Signal Process},
         title={Deep Equilibrium Learning of Explicit Regularizers 
         for Imaging Inverse Problems},
         year={2023},
         note={DOI:10.1109/OJSP.2023.3296036}
        }

@article{Daniele2026,
author = {Daniele, Christian and Villa, Silvia and Vaiter, Samuel and Calatroni, Luca},
title = {Deep Equilibrium Models for Poisson Imaging Inverse Problems via Mirror Descent},
journal = {SIAM Journal on Imaging Sciences},
volume = {19},
number = {2},
pages = {1077-1109},
year = {2026},
doi = {10.1137/25M1779048},
URL = { https://doi.org/10.1137/25M1779048},
eprint = {https://doi.org/10.1137/25M1779048}
}

@misc{learning2026papafitsoros,
      title={Learning spatially varying regularisation parameters of low regularity for image reconstruction}, 
      author={Kostas Papafitsoros and Luca Calatroni and Andreas Kofler},
      year={2026},
      eprint={2608.25127},
      archivePrefix={arXiv},
      primaryClass={eess.IV},
      url={https://arxiv.org/abs/2608.25127}, 
}

@article{pruessmann2001advances,
  author    = {Pruessmann, Klaas P. and Weiger, Markus and B{\"o}rnert, Peter and Boesiger, Peter},
  title     = {Advances in sensitivity encoding with arbitrary k-space trajectories},
  number    = {4},
  pages     = {638--651},
  volume    = {46},
  journal   = {Magnetic Resonance in Medicine},
  publisher = {Wiley Online Library},
  year      = {2001},
}

@article{zhang2017beyond,
  title={Beyond a {G}aussian denoiser: Residual learning of deep {CNN} for image denoising},
  author={Zhang, Kai and Zuo, Wangmeng and Chen, Yunjin and Meng, Deyu and Zhang, Lei},
  journal={IEEE Transactions on Image Processing},
  volume={26},
  number={7},
  pages={3142--3155},
  year={2017},
  publisher={IEEE}
}

@inproceedings{venkatakrishnan2013plug,
  title={Plug-and-play priors for model based reconstruction},
  author={Venkatakrishnan, Singanallur V and Bouman, Charles A and Wohlberg, Brendt},
  booktitle={2013 IEEE Global Conference on Signal and Information Processing},
  pages={945--948},
  year={2013},
  organization={IEEE}
}

@data{mdcnn_data,
author = {El-Rewaidy, Hossam},
publisher = {Harvard Dataverse},
title = {{Replication Data for: Multi-Domain Convolutional Neural Network (MD-CNN) For Radial Reconstruction of Dynamic Cardiac MRI}},
year = {2020},
version = {V2},
doi = {10.7910/DVN/CI3WB6},
url = {https://doi.org/10.7910/DVN/CI3WB6}
}

@ARTICLE{lefkiammiatis2013,
  author={Lefkimmiatis, Stamatios and Ward, John Paul and Unser, Michael},
  journal={IEEE Transactions on Image Processing}, 
  title={Hessian Schatten-Norm Regularization for Linear Inverse Problems}, 
  year={2013},
  volume={22},
  number={5},
  pages={1873-1888},
  doi={10.1109/TIP.2013.2237919}}

@Inbook{kofler2024quantitative,
author="Kofler, Andreas
and Zimmermann, Felix Frederik
and Papafitsoros, Kostas",
editor="Sack, Ingolf
and Schaeffter, Tobias",
title="Machine Learning for Quantitative Magnetic Resonance Image Reconstruction",
bookTitle="Quantification of Biophysical Parameters in Medical Imaging",
year="2024",
publisher="Springer International Publishing",
address="Cham",
pages="171--213",
isbn="978-3-031-61846-8",
doi="10.1007/978-3-031-61846-8_9"
}

@article{zimmermann2025mrpro,
  title={{MRpro-PyTorch-based MR} image reconstruction and processing package},
  author={Zimmermann, Felix Frederik and Schuenke, Patrick and Brahma, Sherine and Guastini, Mara and Hammacher, Johannes and Kofler, Andreas and Kranich Redshaw, Catarina and Lunin, Leonid and Martin, Stefan and Schote, David and others},
  journal={Zenodo},
  year={2025}
}

@article{gilton2021deep,
  title={Deep equilibrium architectures for inverse problems in imaging},
  author={Gilton, Davis and Ongie, Gregory and Willett, Rebecca},
  journal={IEEE Transactions on Computational Imaging},
  volume={7},
  pages={1123--1133},
  year={2021},
  publisher={IEEE}
}

@inproceedings{gregor2010learning,
  title={Learning fast approximations of sparse coding},
  author={Gregor, Karol and LeCun, Yann},
  booktitle={Proceedings of the 27th international conference on international conference on machine learning},
  pages={399--406},
  year={2010}
}

@inproceedings{zhang2018ista,
  title={ISTA-Net: Interpretable optimization-inspired deep network for image compressive sensing},
  author={Zhang, Jian and Ghanem, Bernard},
  booktitle={2018 IEEE/CVF Conference on Computer Vision and Pattern Recognition},
  pages={1828--1837},
  year={2018},
  organization={IEEE}
}

@article{bai2019deep,
  title={Deep equilibrium models},
  author={Bai, Shaojie and Kolter, J Zico and Koltun, Vladlen},
  journal={Advances in Neural Information Processing Systems},
  volume={32},
  year={2019}
}

@article{sun2016deep,
  title={Deep {ADMM-Net} for compressive sensing {MRI}},
  author={Sun, Jian and Li, Huibin and Xu, Zongben and others},
  journal={Advances in Neural Information Processing Systems},
  volume={29},
  year={2016}
}

@article{adler2018learned,
  title={Learned primal-dual reconstruction},
  author={Adler, Jonas and {\"O}ktem, Ozan},
  journal={IEEE Transactions on Medical Imaging},
  volume={37},
  number={6},
  pages={1322--1332},
  year={2018},
  publisher={IEEE}
}

@article{zimmermann2026mrpro_arxiv,
      title={{MRpro}: open framework for model-based, learned, and quantitative {MR} imaging}, 
      author={Felix Frederik Zimmermann and Patrick Schuenke and Christoph S. Aigner and Bill A. Bernhardt and Mara Guastini and Johannes Hammacher and Noah Jaitner and Andreas Kofler and Leonid Lunin and Stefan Martin and Catarina Redshaw Kranich and Jakob Schattenfroh and David Schote and Yanglei Wu and Christoph Kolbitsch},
      year={2026},
      eprint={2507.23129},
      archivePrefix={arXiv},
      primaryClass={eess.IV},
      url={https://arxiv.org/abs/2507.23129}, 
}

@article{lyu2025state,
  title={The state-of-the-art in cardiac {MRI} reconstruction: Results of the {CMRxRecon} challenge in {MICCAI} 2023},
  author={Lyu, Jun and Qin, Chen and Wang, Shuo and Wang, Fanwen and Li, Yan and Wang, Zi and Guo, Kunyuan and Ouyang, Cheng and T{\"a}nzer, Michael and Liu, Meng and others},
  journal={Medical Image Analysis},
  volume={101},
  pages={103485},
  year={2025},
  publisher={Elsevier}
}

@article{sidky2022report,
  title={Report on the {AAPM} deep-learning sparse-view {CT} grand challenge},
  author={Sidky, Emil Y and Pan, Xiaochuan},
  journal={Medical Physics},
  volume={49},
  number={8},
  pages={4935--4943},
  year={2022},
  publisher={Wiley Online Library}
}

@article{muckley2021results,
  title={Results of the 2020 {fastMRI} challenge for machine learning {MR} image reconstruction},
  author={Muckley, Matthew J and Riemenschneider, Bruno and Radmanesh, Alireza and Kim, Sunwoo and Jeong, Geunu and Ko, Jingyu and Jun, Yohan and Shin, Hyungseob and Hwang, Dosik and Mostapha, Mahmoud and others},
  journal={IEEE Transactions on Medical Imaging},
  volume={40},
  number={9},
  pages={2306--2317},
  year={2021},
  publisher={IEEE}
}

@article{beauferris2022multi,
  title={Multi-coil mri reconstruction challenge—assessing brain {MRI} reconstruction models and their generalizability to varying coil configurations},
  author={Beauferris, Youssef and Teuwen, Jonas and Karkalousos, Dimitrios and Moriakov, Nikita and Caan, Matthan and Yiasemis, George and Rodrigues, L{\'\i}via and Lopes, Alexandre and Pedrini, Helio and Rittner, Let{\'\i}cia and others},
  journal={Frontiers in Neuroscience},
  volume={16},
  pages={919186},
  year={2022},
  publisher={Frontiers Media SA}
}

@inproceedings{agustsson2017ntire,
  title={{NTIRE} 2017 challenge on single image super-resolution: Dataset and study},
  author={Agustsson, Eirikur and Timofte, Radu},
  booktitle={2017 IEEE Conference on Computer Vision and Pattern Recognition Workshops (CVPRW)},
  pages={1122--1131},
  year={2017},
  organization={IEEE}
}

@InProceedings{noise2self,
  title = 	 {{N}oise2{S}elf: Blind Denoising by Self-Supervision},
  author =       {Batson, Joshua and Royer, Loic},
  booktitle = 	 {Proceedings of the 36th International Conference on Machine Learning},
  pages = 	 {524--533},
  year = 	 {2019},
  editor = 	 {Chaudhuri, Kamalika and Salakhutdinov, Ruslan},
  volume = 	 {97},
  series = 	 {Proceedings of Machine Learning Research},
  month = 	 {09--15 Jun},
  publisher =    {PMLR}
}

@misc{wang_benchmarking_2026,
	title = {Benchmarking {Self}-{Supervised} {Learning} {Methods} for {Accelerated} {MRI} {Reconstruction}},
	url = {http://arxiv.org/abs/2502.14009},
	doi = {10.48550/arXiv.2502.14009},
	urldate = {2026-03-22},
	publisher = {arXiv},
	author = {Wang, Andrew and McDonagh, Steven and Davies, Mike},
	year = {2026},
	note = {version: 5},
}

@article{li_self-supervised_2025,
	title = {Self-supervised learning for {MRI} reconstruction: a review and new perspective},
	volume = {38},
	issn = {1352-8661},
	shorttitle = {Self-supervised learning for {MRI} reconstruction},
	doi = {10.1007/s10334-025-01274-y},
	language = {en},
	number = {6},
	urldate = {2026-03-22},
	journal = {Magnetic Resonance Materials in Physics, Biology and Medicine},
	author = {Li, Xinzhen and Huang, Jinhong and Sun, Guanglong and Yang, Zihan},
	month = dec,
	year = {2025},
	pages = {1053--1074},
}

@article{elad_image_2023,
	title = {Image Denoising: {The} Deep Learning Revolution and Beyond -- {A} Survey Paper},
	copyright = {arXiv.org perpetual, non-exclusive license},
	shorttitle = {Image {Denoising}},
	url = {https://arxiv.org/abs/2301.03362},
	doi = {10.48550/ARXIV.2301.03362},
	urldate = {2026-02-25},
	publisher = {arXiv},
	author = {Elad, Michael and Kawar, Bahjat and Vaksman, Gregory},
	year = {2023},
	note = {Version Number: 1},
}

@article{zhang_unleashing_nodate,
	title = {Unleashing the Power of Self-Supervised Image Denoising: A Comprehensive Review},
    year={2023},
	language = {en},
    url = {https://arxiv.org/abs/2308.00247},
	author = {Zhang, Dan and Zhou, Fangfang and Albu, Felix and Wei, Yuanzhou and Yang, Xiao and Gu, Yuan and Li, Qiang},
}

@article{pragliola2023and,
  title={On and beyond total variation regularization in imaging: the role of space variance},
  author={Pragliola, Monica and Calatroni, Luca and Lanza, Alessandro and Sgallari, Fiorella},
  journal={SIAM Review},
  volume={65},
  number={3},
  pages={601--685},
  year={2023},
  publisher={SIAM}
}

@Article{holler2014infimal,
  author    = {Holler, Martin and Kunisch, Karl},
  journal   = {SIAM Journal on Imaging Sciences},
  title     = {On infimal convolution of {TV}-type functionals and applications to video and image reconstruction},
  year      = {2014},
  number    = {4},
  pages     = {2258--2300},
  volume    = {7},
  publisher = {SIAM},
}

@article{kofler2023learning,
  title={Learning regularization parameter-maps for variational image reconstruction using deep neural networks and algorithm unrolling},
  author={Kofler, Andreas and Altekr{\"u}ger, Fabian and Antarou Ba, Fatima and Kolbitsch, Christoph and Papoutsellis, Evangelos and Schote, David and Sirotenko, Clemens and Zimmermann, Felix Frederik and Papafitsoros, Kostas},
  journal={SIAM Journal on Imaging Sciences},
  volume={16},
  pages={2202--2246},
  year={2023},
  publisher={SIAM}
}

@article{kolbitsch2014cardiac,
  title={Cardiac functional assessment without electrocardiogram using physiological self-navigation},
  author={Kolbitsch, Christoph and Prieto, Claudia and Schaeffter, Tobias},
  journal={Magnetic Resonance in Medicine},
  volume={71},
  number={3},
  pages={942--954},
  year={2014},
  publisher={Wiley Online Library}
}

@article{kulathilake2023review,
  title={A review on deep learning approaches for low-dose computed tomography restoration},
  author={Kulathilake, KA Saneera Hemantha and Abdullah, Nor Aniza and Sabri, Aznul Qalid Md and Lai, Khin Wee},
  journal={Complex \& Intelligent Systems},
  volume={9},
  number={3},
  pages={2713--2745},
  year={2023},
  publisher={Springer}
}

@article{romano2017little,
  title={The little engine that could: Regularization by denoising ({RED})},
  author={Romano, Yaniv and Elad, Michael and Milanfar, Peyman},
  journal={SIAM Journal on Imaging Sciences},
  volume={10},
  number={4},
  pages={1804--1844},
  year={2017},
  publisher={SIAM}
}

@INPROCEEDINGS{schulz2026learning,
  author={Schulz, Joshua and Schote, David and Kolbitsch, Christoph and Papafitsoros, Kostas and Kofler, Andreas},
  booktitle={2026 IEEE International Conference on Image Processing (ICIP)}, 
  title={Learning Spatially Adaptive Sparsity Level Maps for Arbitrary Convolutional Dictionaries}, 
  year={2026},
  volume={},
  number={},
  pages={1-6}}

@article{reehorst2018regularization,
  title={Regularization by denoising: Clarifications and new interpretations},
  author={Reehorst, Edward T and Schniter, Philip},
  journal={IEEE Transactions on Computational Imaging},
  volume={5},
  number={1},
  pages={52--67},
  year={2018},
  publisher={IEEE}
}

@article{chen2022ai,
  title={{AI}-based reconstruction for fast {MRI}—a systematic review and meta-analysis},
  author={Chen, Yutong and Sch{\"o}nlieb, Carola-Bibiane and Li{\`o}, Pietro and Leiner, Tim and Dragotti, Pier Luigi and Wang, Ge and Rueckert, Daniel and Firmin, David and Yang, Guang},
  journal={Proceedings of the IEEE},
  volume={110},
  number={2},
  pages={224--245},
  year={2022},
  publisher={IEEE}
}

@article{ongie2020deep,
  title={Deep learning techniques for inverse problems in imaging},
  author={Ongie, Gregory and Jalal, Ajil and Metzler, Christopher A and Baraniuk, Richard G and Dimakis, Alexandros G and Willett, Rebecca},
  journal={IEEE Journal on Selected Areas in Information Theory},
  volume={1},
  number={1},
  pages={39--56},
  year={2020},
  publisher={IEEE}
}

@article{lustig2007sparse,
  title={Sparse {MRI}: The application of compressed sensing for rapid {MR} imaging},
  author={Lustig, Michael and Donoho, David and Pauly, John M},
  journal={Magnetic Resonance in Medicine},
  volume={58},
  number={6},
  pages={1182--1195},
  year={2007},
  publisher={Wiley Online Library}
}

@article{sidky2012convex,
  title={Convex optimization problem prototyping for image reconstruction in computed tomography with the {Chambolle}--{Pock} algorithm},
  author={Sidky, Emil Y and J{\o}rgensen, Jakob H and Pan, Xiaochuan},
  journal={Physics in Medicine and Biology},
  volume={57},
  number={10},
  pages={3065--3091},
  year={2012},
  publisher={IOP Publishing}
}

@inproceedings{krull2019noise2void,
  title={Noise2void-learning denoising from single noisy images},
  author={Krull, Alexander and Buchholz, Tim-Oliver and Jug, Florian},
  booktitle={Proceedings of the IEEE/CVF Conference on Computer Vision and Pattern Recognition},
  pages={2129--2137},
  year={2019}
}

@inproceedings{vu2025deep,
  title={Deep unrolling for learning optimal spatially varying regularisation parameters for total generalised variation},
  author={Vu, Thanh Trung and Kofler, Andreas and Papafitsoros, Kostas},
  booktitle={International Conference on Scale Space and Variational Methods in Computer Vision},
  pages={282--294},
  year={2025},
  organization={Springer}
}

@article{wang2025cmrxrecon2024,
  title={{CMRxRecon2024}: a multimodality, multiview k-space dataset boosting universal machine learning for accelerated cardiac {MRI}},
  author={Wang, Zi and Wang, Fanwen and Qin, Chen and Lyu, Jun and Ouyang, Cheng and Wang, Shuo and Li, Yan and Yu, Mengyao and Zhang, Haoyu and Guo, Kunyuan and others},
  journal={Radiology: Artificial Intelligence},
  volume={7},
  number={2},
  pages={e240443},
  year={2025},
  publisher={Radiological Society of North America}
}

@article{wang2024cmrxrecon,
  title={{CMRxRecon}: A publicly available k-space dataset and benchmark to advance deep learning for cardiac {MRI}},
  author={Wang, Chengyan and Lyu, Jun and Wang, Shuo and Qin, Chen and Guo, Kunyuan and Zhang, Xinyu and Yu, Xiaotong and Li, Yan and Wang, Fanwen and Jin, Jianhua and others},
  journal={Scientific Data},
  volume={11},
  number={1},
  pages={687},
  year={2024},
  publisher={Nature Publishing Group UK London}
}

@inproceedings{darestani2022test,
  title={Test-time training can close the natural distribution shift performance gap in deep learning based compressed sensing},
  author={Darestani, Mohammad Zalbagi and Liu, Jiayu and Heckel, Reinhard},
  booktitle={International Conference on Machine Learning},
  pages={4754--4776},
  year={2022},
  organization={PMLR}
}

@inproceedings{
heckel2018deep,
title={Deep Decoder: Concise Image Representations from Untrained Non-convolutional Networks},
author={Reinhard Heckel and Paul Hand},
booktitle={International Conference on Learning Representations},
year={2019}
}

@inproceedings{liu2019image,
  title={Image restoration using total variation regularized deep image prior},
  author={Liu, Jiaming and Sun, Yu and Xu, Xiaojian and Kamilov, Ulugbek S},
  booktitle={ICASSP 2019-2019 IEEE International Conference on Acoustics, Speech and Signal Processing (ICASSP)},
  pages={7715--7719},
  year={2019}
}

@article{monga2021algorithm,
  title={Algorithm unrolling: Interpretable, efficient deep learning for signal and image processing},
  author={Monga, Vishal and Li, Yuelong and Eldar, Yonina C},
  journal={IEEE Signal Processing Magazine},
  volume={38},
  number={2},
  pages={18--44},
  year={2021},
  publisher={IEEE}
}

@inproceedings{
yaman2022zeroshot,
title={Zero-Shot Self-Supervised Learning for {MRI} Reconstruction},
author={Burhaneddin Yaman and Seyed Amir Hossein Hosseini and Mehmet Akcakaya},
booktitle={International Conference on Learning Representations},
year={2022}
}

@inproceedings{ulyanov2018deep,
  title={Deep image prior},
  author={Ulyanov, Dmitry and Vedaldi, Andrea and Lempitsky, Victor},
  booktitle={Proceedings of the IEEE conference on computer vision and pattern recognition},
  pages={9446--9454},
  year={2018}
}

@INPROCEEDINGS{kofler2025ell1,
  author={Kofler, Andreas and Calatroni, Luca and Kolbitsch, Christoph and Papafitsoros, Kostas},
  booktitle={2025 33rd European Signal Processing Conference (EUSIPCO)},
  title={Learning Spatially Adaptive $\ell_{1}$-Norms Weights for Convolutional Synthesis Regularization},
  year={2025},
  volume={},
  number={},
  pages={1782-1786},
  doi={10.23919/EUSIPCO63237.2025.11226665}}

@article{ChambolleLions1997,
	Author = {Chambolle, Antonin and Lions, Pierre-Louis},
	Doi = {10.1007/s002110050258},
	Isbn = {0945-3245},
	Journal = {Numerische Mathematik},
	Number = {2},
	Pages = {167--188},
	Title = {Image recovery via total variation minimization and related problems},
	Ty = {JOUR},
	Url = {https://doi.org/10.1007/s002110050258},
	Volume = {76},
	Year = {1997}}

@article{ChambollePock2011FirstOrderPrimalDual,
  title   = {A First-Order Primal-Dual Algorithm for Convex Problems with Applications to Imaging},
  author  = {Chambolle, Antonin and Pock, Thomas},
  journal = {Journal of Mathematical Imaging and Vision},
  volume  = {40},
  number  = {1},
  pages   = {120--145},
  year    = {2011},
  doi     = {10.1007/s10851-010-0251-1}
}

@article{BrediesKunischPock2010TGV,
  title   = {{Total Generalized Variation}},
  author  = {Bredies, Kristian and Kunisch, Karl and Pock, Thomas},
  journal = {SIAM Journal on Imaging Sciences},
  volume  = {3},
  number  = {3},
  pages   = {492--526},
  year    = {2010},
  doi     = {10.1137/090769521}
}

@article{yaman2020self,
  title={Self-supervised learning of physics-guided reconstruction neural networks without fully sampled reference data},
  author={Yaman, Burhaneddin and Hosseini, Seyed Amir Hossein and Moeller, Steen and Ellermann, Jutta and U{\u{g}}urbil, K{\^a}mil and Ak{\c{c}}akaya, Mehmet},
  journal={Magnetic Resonance in Medicine},
  volume={84},
  number={6},
  pages={3172--3191},
  year={2020},
  publisher={Wiley Online Library}
}

@article{ronneberger2015u,
  author       = {Ronneberger, Olaf and Fischer, Philipp and Brox, Thomas},
  title        = {U-net: Convolutional networks for biomedical image segmentation},
  pages        = {234--241},
  volume       = {part III 18},
  journal      = {Medical Image Computing and Computer-Assisted Intervention (MICCAI)},
  organization = {Springer},
  year         = {2015},
}

@Article{mdcnn,
  author  = {El-Rewaidy, Hossam and Fahmy, Ahmed S. and Pashakhanloo, Farhad and Cai, Xiaoying and Kucukseymen, Selcuk and Csecs, Ibolya and Neisius, Ulf and Haji-Valizadeh, Hassan and Menze, Bjoern and Nezafat, Reza},
  journal = {Magnetic Resonance in Medicine},
  title   = {Multi-domain convolutional neural network ({MD-CNN}) for radial reconstruction of dynamic cardiac {MRI}},
  number  = {3},
  pages   = {1195--1208},
  volume  = {85},
  doi     = {10/g9vfqt},
}

@Article{walsh_adaptive_2000,
  author  = {Walsh, David O. and Gmitro, Arthur F. and Marcellin, Michael W.},
  journal = {Magnetic Resonance in Medicine},
  title   = {Adaptive reconstruction of phased array {MR} imagery},
  year    = {2000},
  issn    = {07403194},
  number  = {5},
  pages   = {682--690},
  volume  = {43},
  doi     = {10/bshdwg},
  pmid    = {10800033},
}

@Article{uecker2013,
  author    = {Uecker, Martin and Lai, Peng and Murphy, Mark J. and Virtue, Patrick and Elad, Michael and Pauly, John M. and Vasanawala, Shreyas S. and Lustig, Michael},
  journal   = {Magnetic Resonance Imaging},
  title     = {{ESPIRiT}-an Eigenvalue Approach to Autocalibrating Parallel {MRI}: Where {SENSE} Meets {GRAPPA}.},
  year      = {2013},
  doi       = {10/gfvjn3},
  pmid      = {23649942},
  publisher = {Wiley},
}

@article{tachella2026self,
  title   = {Self-Supervised Learning from Noisy and Incomplete Data},
  author  = {Tachella, Juli{\'a}n and Davies, Mike},
  journal = {Foundations and Trends in Signal Processing},
  year    = {2026},
  volume  = {20},
  number  = {2},
  pages   = {85--184},
  doi     = {10.1108/FTSIG-10-2025-0133},
  url     = {https://doi.org/10.1108/FTSIG-10-2025-0133}
}

@article{aggarwal2018modl,
  author    = {Aggarwal, Hemant K. and Mani, Merry P. and Jacob, Mathews},
  title     = {{MoDL}: Model-based deep learning architecture for inverse problems},
  number    = {2},
  pages     = {394--405},
  volume    = {38},
  journal   = {IEEE Transactions on Medical Imaging},
  publisher = {IEEE},
  year      = {2018},
}

@Article{wang2004image,
  author    = {Wang, Zhou and Bovik, Alan Conrad and Sheikh, Hamid Rahim and Simoncelli, Eero P.},
  journal   = {IEEE Transactions on Image Processing},
  title     = {{Image quality assessment: From error visibility to structural similarity}},
  year      = {2004},
  issn      = {10577149},
  number    = {4},
  pages     = {600--612},
  volume    = {13},
  doi       = {10/c7sr27},
  pmid      = {15376593},
  publisher = {IEEE},
}

@article{tachella2025deepinverse,
    title = {{DeepInverse}: A {Python} package for solving imaging inverse problems with deep learning},
    journal = {Journal of Open Source Software},
    doi = {10.21105/joss.08923},
    url = {https://doi.org/10.21105/joss.08923},
    year = {2025},
    publisher = {The Open Journal},
    volume = {10},
    number = {115},
    pages = {8923},
    author = {Tachella, Julián and Terris, Matthieu and Hurault, Samuel and Wang, Andrew and Davy, Leo and Scanvic, Jérémy and Sechaud, Victor and Vo, Romain and Moreau, Thomas and Davies, Thomas and Chen, Dongdong and Laurent, Nils and Monroy, Brayan and Dong, Jonathan and Hu, Zhiyuan and Nguyen, Minh-Hai and Sarron, Florian and Weiss, Pierre and Escande, Paul and Massias, Mathurin and Modrzyk, Thibaut and Levac, Brett and Liaudat, Tobías I. and Song, Maxime and Hertrich, Johannes and Neumayer, Sebastian and Schramm, Georg},
}

@INPROCEEDINGS{Roth2005,
  author={Roth, S. and Black, M.J.},
  booktitle={2005 IEEE Computer Society Conference on Computer Vision and Pattern Recognition (CVPR'05)}, 
  title={Fields of Experts: a framework for learning image priors}, 
  year={2005},
  volume={2},
  number={},
  pages={860-867 vol. 2},
  doi={10.1109/CVPR.2005.160}}

@inproceedings{pourya2025dealing,
title={{DEAL}ing with Image Reconstruction: Deep Attentive Least Squares},
author={Mehrsa Pourya and Erich Kobler and Michael Unser and Sebastian Neumayer},
booktitle={Forty-second International Conference on Machine Learning},
year={2025},
url={https://openreview.net/forum?id=mMasOShOVt}
}

@article{reisenhofer2018haar,
  title={A {Haar} wavelet-based perceptual similarity index for image quality assessment},
  author={Reisenhofer, Rafael and Bosse, Sebastian and Kutyniok, Gitta and Wiegand, Thomas},
  journal={Signal Processing: Image Communication},
  volume={61},
  pages={33--43},
  year={2018},
  publisher={Elsevier}
}

@inproceedings{hosseini2020high,
  title={High-fidelity accelerated {MRI} reconstruction by scan-specific fine-tuning of physics-based neural networks},
  author={Hosseini, Seyed Amir Hossein and Yaman, Burhaneddin and Moeller, Steen and Ak{\c{c}}akaya, Mehmet},
  booktitle={2020 42nd Annual International Conference of the IEEE Engineering in Medicine \& Biology Society (EMBC)},
  pages={1481--1484},
  year={2020},
  organization={IEEE}
}

@inproceedings{
loshchilov2018decoupled,
title={Decoupled Weight Decay Regularization},
author={Ilya Loshchilov and Frank Hutter},
booktitle={International Conference on Learning Representations},
year={2019}
}

@inproceedings{
terris2026ram,
title={{Reconstruct Anything Model} a lightweight general model for computational imaging},
author={Matthieu Terris and Samuel Hurault and Maxime Song and Juli{\'a}n Tachella},
booktitle={The Fourteenth International Conference on Learning Representations},
year={2026},
url={https://openreview.net/forum?id=Ks9zNS6OsU}
}

@article{crockett2022bilevel,
  title={Bilevel methods for image reconstruction},
  author={Crockett, Caroline and Fessler, Jeffrey A},
  journal={Foundations and Trends{\textregistered} in Signal Processing},
  volume={15},
  number={2-3},
  pages={121--289},
  year={2022},
  publisher={Emerald Publishing Limited}
}

@InProceedings{kingmaAdamMethodStochastic2015a,
  author     = {Kingma, Diederik P. and Ba, Jimmy},
  booktitle  = {3rd {{International Conference}} on {{Learning Representations}}, {{ICLR}} 2015, {{San Diego}}, {{CA}}, {{USA}}, {{May}} 7-9, 2015, {{Conference Track Proceedings}}},
  title      = {Adam: {{A Method}} for {{Stochastic Optimization}}},
  year       = {2015},
  editor     = {Bengio, Yoshua and LeCun, Yann},
  eprint     = {1412.6980},
  eprinttype = {arxiv},
  shorttitle = {Adam},
}

@article{chen2016memory,
  author       = {Tianqi Chen and
                  Bing Xu and
                  Chiyuan Zhang and
                  Carlos Guestrin},
  title        = {Training Deep Nets with Sublinear Memory Cost},
  journal      = {CoRR},
  volume       = {abs/1604.06174},
  year         = {2016},
  url          = {http://arxiv.org/abs/1604.06174},
  eprinttype   = {arXiv},
  eprint       = {1604.06174},
  bibsource    = {dblp computer science bibliography, https://dblp.org}
}
\bibliographystyle{IEEEtran}
\end{document}